\documentstyle[11pt]{article}

\input xy
\xyoption{all}

\font\tengoth=eufb10
\font\sevengoth=eufb7
\font\fivegoth=eufb5
\newfam\gothfam
\textfont\gothfam=\tengoth
\scriptfont\gothfam=\sevengoth
\scriptscriptfont\gothfam=\fivegoth
\def\goth{\fam\gothfam\tengoth}

\def\C{{\mathchoice {\setbox0=\hbox{$\displaystyle\mathrm C$}\hbox{\hbox
to0pt{\kern0.4\wd0\vrule height0.9\ht0\hss}\box0}}
{\setbox0=\hbox{$\textstyle\mathrm C$}\hbox{\hbox
to0pt{\kern0.4\wd0\vrule height0.9\ht0\hss}\box0}}
{\setbox0=\hbox{$\scriptstyle\mathrm C$}\hbox{\hbox
to0pt{\kern0.4\wd0\vrule height0.9\ht0\hss}\box0}}
{\setbox0=\hbox{$\scriptscriptstyle\mathrm C$}\hbox{\hbox
to0pt{\kern0.4\wd0\vrule height0.9\ht0\hss}\box0}}}}
\def\R{{\mathrm{I\!R}}} 
\def\N{{\mathrm{I\!N}}} 

\def\Q{{\mathchoice {\setbox0=\hbox{$\displaystyle\mathrm
Q$}\hbox{\raise
0.15\ht0\hbox to0pt{\kern0.4\wd0\vrule height0.8\ht0\hss}\box0}}
{\setbox0=\hbox{$\textstyle\mathrm Q$}\hbox{\raise
0.15\ht0\hbox to0pt{\kern0.4\wd0\vrule height0.8\ht0\hss}\box0}}
{\setbox0=\hbox{$\scriptstyle\mathrm Q$}\hbox{\raise
0.15\ht0\hbox to0pt{\kern0.4\wd0\vrule height0.7\ht0\hss}\box0}}
{\setbox0=\hbox{$\scriptscriptstyle\mathrm Q$}\hbox{\raise
0.15\ht0\hbox to0pt{\kern0.4\wd0\vrule height0.7\ht0\hss}\box0}}}}
\def\Z{{\mathchoice {\hbox{$\sf\textstyle Z\kern-0.4em Z$}}
{\hbox{$\sf\textstyle Z\kern-0.4em Z$}}
{\hbox{$\sf\scriptstyle Z\kern-0.3em Z$}}
{\hbox{$\sf\scriptscriptstyle Z\kern-0.2em Z$}}}}

\newcounter{amoi}
\newtheorem{theorem}{Theorem}[subsection]
\newtheorem{lemma}[theorem]{Lemma}
\newtheorem{proposition}[theorem]{Proposition}
\newtheorem{corollary}[theorem]{Corollary}
\newtheorem{definition}[theorem]{Definition}

\def\hfl#1#2{\smash{\mathop{\hbox to 12 mm{\rightarrowfill}}
\limits^{\scriptstyle#1}_{\scriptstyle#2}}}

\begin{document}

\title{Formality theorem with  coefficients in a  module}
\author{Sophie Chemla}

\maketitle

\begin{abstract}
In this article, $X$ will denote a ${\cal C}^{\infty}$ manifold. 
In a very famous article, Kontsevich (~\cite{Ko}) 
showed that the differential
graded Lie algebra (DGLA) of polydifferential operators on $X$ is formal.
Calaque (~\cite{C1}) extended this theorem to any Lie algebroid. More
precisely, given any  Lie algebroid $E$ over $X$, he  defined the DGLA
of $E$-polydifferential operators, 
$\Gamma \left (X,  ^{E}D^{*}_{poly} \right )$, 
and showed that it is formal. 
Denote by $\Gamma \left ( X, ^{E}T^{*}_{poly}\right )$ 
the DGLA of $E$-polyvector fields. 
Considering $M$, a module over $E$, we define 
$\Gamma \left (X,  ^{E}T_{poly}^{*}(M)\right )$ the 
$\Gamma \left ( X, ^{E}T^{*}_{poly}\right )$-module of $E$-polyvector
fields with values in $M$. 
Similarly, we define the 
$\Gamma \left ( X, ^{E}D^{*}_{poly}\right )$-module 
of $E$-polydifferential operators 
with values in $M$, $\Gamma \left (X,  ^{E}D^{*}_{poly}(M)\right )$. 
We show that there is a 
quasi-isomorphism of $L_{\infty}$-modules over 
$\Gamma \left ( X, ^{E}T^{*}_{poly}\right )$ from 
$\Gamma \left ( X, ^{E}T^{*}_{poly}(M)\right )$ to 
$\Gamma \left ( X, ^{E}D^{*}_{poly}(M)\right )$. 
Our result extends  Calaque 's 
(and Kontsevich's) result.  
\end{abstract}

\section{Introduction}

In this article, $X$ will denote a ${\cal C}^{\infty}$-manifold and
${\cal O}_{X}$ will denote the sheaf of ${\cal C}^{\infty}$
functions. 
To $X$ are associated two sheaves of 
differential graded Lie algebras (DGLAs ) 
$T^{*}_{poly}$ and $D^{*}_{poly}$. The first one, $T^{*}_{poly}$ is
the sheaf of DGLAs of polyvector fields on $X$ with differential zero and
Schouten bracket. The second one, $D^{*}_{poly}$, is the sheaf of DGLAs of
polydifferential operators on $X$ with Hochschild differential and 
Gerstenhaber bracket. Kontsevich showed that there is a
quasi-isomorphism of $L_{\infty}$-algebras from 
$\Gamma \left ( X, T^{*}_{poly} \right )$ to
$\Gamma \left (X,  D^{*}_{poly} \right )$, 
that is to say that $\Gamma (X,D^{*}_{poly})$ is formal. The aim of
this article is to introduce a module in Kontsevitch formality
theorem. 
 
Let us now consider a ${\cal D}_{X}$-module $M$. Inspired by the
expression of the Schouten bracket, we endow 
$T^{*}_{poly}(M)=
T^{*}_{poly}{\displaystyle \mathop \otimes_{{\cal O}_{X}} M}$ with a 
$T^{*}_{poly}$-module structure. 
Similarly, we can endow 
$D^{*}_{poly}(M)=D^{*}_{poly}
{\displaystyle \mathop \otimes _{{\cal O}_{X}}} M$ with a 
$D_{poly}^{*}$-module structure as follows : 
if $P \in D_{poly}^{p}$ and $Q \in D_{poly}^{q}(M)$, 
$$P \cdot_{G} Q = P \bullet Q -
(-1)^{pq} Q \bullet P $$
 with 
$$\begin{array}{l}
\forall a_{0}, \dots, a_{p+q} \in {\cal O}_{X},\\
(P \bullet Q ) (a_{0}, \dots, a_{p+q})=
{\displaystyle \sum _{i=0}^{p}}(-1)^{iq}
P(a_{0}, \dots, a_{i-1}, Q(a_{i}, \dots, a_{i+q}), \dots, a_{p+q}).
\end{array}$$
The formula makes sense because $Q$ is a differential operator with
coefficients in a ${\cal D}_{X}$-module $M$. The expression 
$Q \bullet P$ is defined in an analogous way. 
The differential on $D_{poly}^{*}(M)$ is given by the action of the
multiplication $\mu$, $\mu \cdot_{G} -$. 
Using Kontsevich's formality theorem, one may see $D^{*}_{poly}(M)$
 as an $L_{\infty}$-module over $T^{*}_{poly}$ and we will prove that it
 is formal. We will work in the more general setting of Lie algebroids.

Let us now consider a Lie algebroid $E$. To $E$ is associated a sheaf
of $E$-differential operators,  $D(E)$ (~\cite{R}). 
Lie algebroids generalize at the same time  the sheaf of vector fields
on a manifold (in this case $E=TX$ and $D(E)={\cal D}_{X}$) and 
Lie algebras (in this case $D(E)$ is the enveloping algebra). 
Lie algebroids have been
extensively studied recently because many examples of Lie algebroids
arise from geometry (Poisson manifolds, group actions, foliations
...).
To $E$, one can associate the sheaf of DGLAs of $E$-polyvectorfields 
$^{E}T^{*}_{poly}= 
{\displaystyle \mathop \oplus_{k=-1}^{\infty}}\wedge ^{k+1}E$
with zero differential and a Schouten type Lie bracket (~\cite{C1}). 
Calaque  has given an appropriate generalization of the notion
of polydifferential operators. In ~\cite{C1}, he defines the DGLA of
$E$-polydifferential operators , 
$\Gamma\left (X, ^{E}D_{poly}^{*}\right )$, 
 and constructs an $L_{\infty}$ quasi-isomorphism from 
$\Gamma \left ( X, ^{E}T_{poly}^{*}\right )$ to 
$\Gamma\left (X, ^{E}D_{poly}^{*}\right )$. 

Let us now consider a $D(E)$-module $M$. We can perform the
construction described above and define 
the $^{E}T_{poly}^{*}$-module $^{E}T_{poly}^{*}(M)$ (the sheaf
of the 
$E$-polyvectors with coefficients in $M$) and 
the $^{E}D^{*}_{poly}$-module $^{E}D_{poly}^{*}(M)$ 
(the sheaf of $E$-polydifferential operators with coefficients in
$M$).
By Calaque 's result  $\Gamma\left (X, ^{E}D_{poly}^{*}(M)\right )$ is an 
$L_{\infty}$-module over  $\Gamma \left ( X, ^{E}T_{poly}^{*}\right )$. 
The main result of the paper is the following theorem.\\

{\bf Theorem 3.4.1 :}

{\it There is a quasi-isomorphism of $L_{\infty}$-modules over 
$\Gamma \left (X, ^{E}T_{poly}\right )$ 
from \linebreak 
$\Gamma \left ( X, ^{E}T_{poly}(M)\right )$ to 
$\Gamma \left (X, ^{E}D_{poly}(M)\right ).$}\\

Our result extends Calaque 's formality theorem (~\cite{C1},  
take $M={\cal O}_{X}$) and Kontsevich's formality theorem 
(~\cite{Ko}, take  
$M={\cal O}_{X}$ and $E=TX$).\\

If $X$ is a Poisson manifold, we know from Kontsevitch's work 
(~\cite{Ko})
that there is a star product on $O=\Gamma ({\cal O}_{X})$. 
Let ${\cal  M}$ be a
${\cal D}_{X}$-module and $M=\Gamma \left (X, {\cal M}\right )$. 
Using the star product, we can endow $M[[h]]$ with a 
$O[[h]]\otimes O[[h]]^{op}$-module structure. 
If $\pi \in \Gamma (\bigwedge ^{2}TX)$ is the bivector defining the
Poisson structure on $O$, $h \pi$ defines a Poisson structure on the
algebra $O[[h]]$. 
As a
corollary of our theorem, we get an isomorphism from the Poisson
cohomology of the Poisson algebra $O[[h]]$ with
coefficients in $M[[h]]$ and the differential Hochschild cohomology of 
$O[[h]]$ with coefficients in $M[[h]]$.\\

Our proofs are analogous to that of 
~\cite{D1},~\cite{C1}, ~\cite{D2}, ~\cite{CDH}. We use 
Kontsevitch's formality theorem for $\R^{d}_{formal}$
and a Fedosov like globalization techniques.\\

{\bf Acknowledgements :} 

I am grateful to D. Calaque, M. Duflo,  
B. Keller, P. Schapira and C. Torossian for helpful discussions. 
I thank D. Calaque and V. Dolgushev for making
comments on  this article. \\

{\bf Notation :}

For a study of $L_{\infty}$ structures, 
we refer to ~\cite{AMM}, \cite{D2}, ~\cite{D3}, 
~\cite{HS}, ~\cite{LS}.\\

Let $k$ be a field of characteristic zero and let $V$ be a $\Z$-graded 
$k$-vector space 
$$V=\mathop \oplus_{i \in \Z}V_{i}.$$
If $x$ is in $V_{i}$, we set $\mid x \mid =i$.
We will always assume that the gradation is bounded below. Recall the
definition of the graded symmetric algebra and the graded wedge
algebra :
$$\begin{array}{l}
S(V)={\displaystyle \frac{T(V)}{<x \otimes y -
(-1)^{\mid x \mid \mid y \mid}y \otimes x >}}\\
\wedge (V)={\displaystyle \frac{T(V)}{<x \otimes y +
(-1)^{\mid x \mid \mid y \mid}y \otimes x >}}.
\end{array}$$
If $i$ is in $\Z$, we will denote by $V[i]$ the graded vector space
defined by $V[i]^{n}=V^{i+n}$.

Denote by $S^{c}(V)$ the cofree cocommutative coalgebra without
counity cofreely
cogenerated by $V$. As a vector space $S^{c}(V)$ is $S^{+}(V)$. Its
comultiplication is given by 
$$\Delta (x_{1}\dots x_{n})=
{\displaystyle \sum_{
\begin{array}{c}
I \sqcup J=[1,n]\\
I \neq \emptyset\\
J\neq \emptyset
\end{array}}}
(-1)^{\epsilon(I ,J)}x_{I} \otimes x_{J}.$$
where $\epsilon(I,J)$ is the number of inversion of odd elements when 
going from $x_{I}x_{J}$ to $x_{1}\dots x_{n}$. 
A coderivation
 $Q$ on $S^{c}(V)$ is determined by its Taylor coefficients 
$Q^{[n]} : S^{n}(V) \to V$ (obtained by composing $Q$ with the
projection from $S(V)$ onto $V$).

An $L_{\infty}$ algebra is a couple $(L, Q)$ where $L$ is a graded
vector space and $Q$ is a degree 1 two-nilpotent coderivation of 
$S^{c}(L[1])=C(L)$. The coderivation $Q$ is determined by its Taylor
coefficients $\left (Q^{[n]}\right )_{n \geq 1}$. Using an isomorphism 
between  $S^{n}(L[1])$ and $\wedge ^{n}(L)[n]$, the Taylor
coefficients may be seen as maps \linebreak
$\overline{Q}^{[n]}: \wedge ^{n}L \to L [2-n]$.
A differential graded Lie algebra $(L,d, [\;,\;])$ (with differential $d$
and Lie bracket $[,]$) gives rise to an 
$L_{\infty}$-algebra determined by $\overline{Q}^{[1]}=d$, 
$\overline{Q}^{[2]}=[,]$ and $\overline{Q}^{[i]}=0$ for $i \geq 2$. 

Let $L$ be a differential graded Lie algebra. 
We will say that it is a filtered DGLA if it is equipped with a complete 
descending filtration $\dots {\cal F}^{1}L \subset {\cal F}^{0}L=L$
such that 
$L={\displaystyle \lim_{n} {L/{\cal F}^{n}L}}$.  
A Maurer Cartan element
of $L$ is an element $x$ of ${\cal F}^{1}L^{1}$ such that 
$Q^{[1]}x + {\displaystyle \frac{1}{2}} Q^{[2]}(x^{2})=0$. 

Let $(L_{1}, Q_{1})$ and $(L_{2}, Q_{2})$ be two
$L_{\infty}$-algebras.
 An $L_{\infty}$-morphism   $F$ from $(L_{1},Q_{1})$ to 
$(L_{2}, Q_{2})$ is a morphism of coalgebras $F: C(L_{1}) \to
C(L_{2})$ compatible with coderivations (this means that 
$F \circ Q_{1}=Q_{2} \circ F$). As $F$ is a morphism of coalgebras, it
is determined by its Taylor coefficients 
$\left ( F^{[n]} : S^{n}(L_{1}[1]) \to L_{2}[1] \right )_{n \geq 1}$ or 
 $\left ( \overline {F}^{[n]} : \wedge^{n}(L_{1}) \to L_{2}[1-n] \right )_{n
   \geq 1}$. The relation $F \circ Q_{1}=Q_{2} \circ F$ boils down 
 to say that $F^{[n]}$ satisfy an infinite collection of equations. 

Let $(L_{1}, Q_{1})$ and $(L_{2}, Q_{2})$ be two filtered DGLAs and 
let $F$ be an $L_{\infty}$-morphism   from $(L_{1}, Q_{1})$
to $(L_{2},Q_{2})$ compatible with these filtrations. 
If $x$ is a Maurer Cartan element of $L_{1}$, then 
${\displaystyle \sum_{n \geq 1} \frac{F^{[n]}(x^{n})}{n!}}$ is a
Maurer Cartan element of $L_{2}$.

Let  $L$ be an $L_{\infty}$-algebra and $M$ a graded vector space. 
We will consider the $C(L)$-comodule $S(L[1])\otimes M$ with the
coaction 
$${\goth a} (x_{1}\dots x_{n}\otimes v)=
{\displaystyle \sum_{
\begin{array}{c}
I \sqcup J=[1,n]\\
I \neq \emptyset\\
\end{array}}} (-1)^{\epsilon (I,J)} x_{I} \otimes (x_{J} \otimes v)$$
where $\epsilon(I,J)$ is the number of inversion of odd elements when 
going from $x_{I}x_{J}$ to $x_{1}\dots x_{n}$.
An $L_{\infty}$-module is a couple 
$(M, \phi )$ where $\phi$ is a degree 1 two-nilpotent 
coderivation of the $C(L)$-comodule  $S(L[1])\otimes M$. 
The coderivation $\phi$ is determined by its Taylor coefficients  
$\phi^{[n]}: S^{n}(L[1]) \otimes M \to M[1]$ or 
$\overline{\phi}^{[n]}: \wedge^{n}(L) \otimes M \to M[1-n]$. 
 The map $\phi ^{[0]}$ is a
differential on $M$. A module $M$ over a differential graded Lie algebra 
$(L,d, [,])$ is an $L_{\infty}$-module with Taylor coefficients 
$\overline{\phi}^{[0]}=d$, $\overline{\phi}^{[1]}(X \otimes m)=X \cdot m$ 
($X \in L$, $m \in M$) and $\overline{\phi}^{[n]}=0$ if $n>1$. 

Let $(M_{1}, \phi_{1})$ and $(M_{2}, \phi_{2})$ be two 
$L_{\infty}$-modules. An $L_{\infty}$-morphism ${\cal V}$
from $(M_{1}, \phi_{1})$ 
to $(M_{2},\phi_{2})$ is a (degree 0) morphism of comodules from 
$S(L[1])\otimes M_{1}$ to $S(L[1]) \otimes M_{2}$ such that 
${\cal V}\circ \phi_{1}=\phi_{2}\circ {\cal V}$. 
It is determined by its Taylor coefficients 
$\left ( {\cal V}^{[n]}: S^{n}(L[1])\otimes M_{1} \to M_{2}
\right )_{n \geq 0}$ or   
 $\left ( {\overline{{\cal V}}^{[n]}}: 
\wedge^{n}(L)\otimes M_{1} \to M_{2}[-n]
\right )_{n \geq 0}$. The compatibility of ${\cal V}$ with
coderivation is expressed by an infinite collection of equations
 satisfied by ${\cal V}^{[n]}$. \\

In this text, DGLA (repectively DGAA) will stand for differential
graded  Lie algebra 
(respectively differential graded associative algebra).   

We assume Einstein convention for the summation over repeated indices.

If ${\cal F}$ is a sheaf over $X$, then $\Gamma ({\cal F})$ denotes
its global sections. If ${\cal F}$ and ${\cal G}$ are two sheaves and
if $\Theta : {\cal F}\to {\cal G}$ is a morphism of sheaves, then 
$\Theta (X)$ will denote the morphism from $\Gamma ({\cal F})$ 
to $\Gamma ({\cal G})$ contained in $\Theta$ 

\section{Recollections}

\subsection{Lie algebroids : definitions and first properties}

Let $X$ be a ${\cal C}^{\infty}$-manifold and let ${\cal O}_{X}$ be the sheaf
of ${\cal C}^{\infty}$ functions on $X$.
Let $\Theta _{X}$ be the ${\cal O}_{X}$-module of ${\cal C}^{\infty}$
 vector fields on $X$.

\begin{definition}
A sheaf in $\R$-Lie algebras  over $X$, $E$,  
is a sheaf of $\R$-vector
spaces such that for any open subset $U$, $E(U)$ is equipped
with the structure of a Lie algebra and  the restriction morphisms are
Lie algebra homomorphisms.
\end{definition}
A morphism between two sheaves of Lie algebras $E$ and
$F$ is a $\R_{X}$-module morphism  which is a Lie algebra
morphism on each open subset.
\begin{definition}
A Lie algebroid over $X$ is a pair
$\left ( E, \omega \right )$ where
\begin{itemize}
\item $E$ is a locally free  ${\cal O}_{X}$-module of finite
 constant rank that is to say a vector bundle over $X$,
\item $E$ is a sheaf of $\R$-Lie algebras,
\item
$\omega : \; E \to \Theta_{X}$  
is an ${\cal O}_{X}$-linear morphism of sheaves  of 
$\R$-Lie algebras  
such that the following
compatibility relation holds :
$$\forall (\xi, \zeta ) \in E^{2}, \;\; \forall f \in
{\cal O}_{X}, \;\;
[ \xi, f \zeta ]= \omega (\xi) (f) \zeta + f [\xi, \zeta].$$

\end{itemize}
\end{definition}
One calls $\omega$  the anchor map. When there is no ambiguity, we
will drop the anchor map in the notation of the Lie algebroid.

For example
 $TX$ is a Lie algebroid on $X$ and a finite dimensional 
Lie algebra is a Lie algebroid over a point. Other examples arises
from Poisson manifolds, foliations, Lie group actions 
(see ~\cite{F} for example).

A Lie algebroid $\left (E , \omega \right )$ gives rise to the
sheaf of $E$-differential operators generated by ${\cal O}_{X}$ and 
$E$ 
which is denoted by 
$D(E)$.
\begin{definition}
$D(E)$ is the sheaf associated to the
presheaf:
$$U \mapsto T^{+}_{\R}
\left ( {\cal O}_{X}(U) \oplus E(U)\right )/J_{U}$$ where
$J_{U}$ is the two sided ideal generated by the relations
$$\begin{array}{l}
\forall (f,g) \in {\cal O}_{X}(U), \;\;
\forall (\xi, \zeta) \in E(U)^{2}\\
\\
1) f \otimes g= fg \\
2) f \otimes \xi =f\xi  \\
3) \xi \otimes \zeta - \zeta \otimes \xi= [ \xi, \zeta ]\\
4) \xi \otimes f - f \otimes \xi = \omega (\xi)(f)
\end{array}$$
\end{definition}

If $E=TX$, $D(E)$ is the sheaf of differential operators on $X$,
${\cal D}_{X}$. If $E$ is a finite dimensional Lie algebra 
$\goth{g}$, 
$D(E)$ is $U(\goth{g})$, 
the enveloping algebra of $\goth{g}.$

$D(E)$ is also endowed with a coassociative ${\cal O}_{X}$-linear coproduct 
$\Delta : D(E) \to D(E) 
{\displaystyle \mathop \otimes_{{\cal O}_{X}}} D(E)$ defined as
follows (see ~\cite{X} example 3.1):
$$\begin{array}{l}
\Delta (1)=1 \otimes 1 \\
\forall u \in E , \;\; \Delta (u)= u \otimes 1 + 1 \otimes u  \\
\forall (P, Q) \in D(E)^{2}, \Delta (PQ)=\Delta (P) \Delta (Q). 
\end{array}$$

Let $M$ be a $D(E)$-module. 
The cohomology of $E$ with coefficients in $M$, 
is computed by the complex  
$\left ( Hom _{{\cal O}_{X}}(\wedge^{*}E, M), {}^{E}d_{M} \right )$ 
where $^{E}d_{M}$ is given by \linebreak  
$\forall \phi \in Hom _{{\cal O}_{X}}\left ( \wedge^{n} E, M \right ), 
\forall u_{0}, \dots,u_{n} \in E$ , 
$$\begin{array}{rcl}
^{E}d_{M} \phi(u_{0}, \dots ,u_{n}) &=&
{\displaystyle \sum_{i=1}^{n}}(-1)^{i}
u_{i}\cdot \phi (u_{1}, \dots, ,\widehat{u_{i}}, \dots, u_{n})\\
&+& {\displaystyle \sum_{i<j}} (-1)^{i+j}
\phi ([u_{i}, u_{j}], u_{0}, \dots, \widehat{u_{i}}, \dots, \widehat{u_{j}},
\dots, u_{n}) 
\end{array}$$

Recall that ${\cal O}_{X}$ has a natural  left $D(E)$-module
structure defined by : 
$$\forall f \in {\cal O}_{X}, \forall P \in D(E), \;\;
P \cdot f =\omega (P)(f).$$
If $M={\cal O}_{X}$, we set $^{E}d_{M}={}^{E}d$ and the complex above
 will be called the Lie  cohomology complex of $E$.\\
  
 If $M$ is a $D(E)$-module, a tensor with coefficients in $M$ is a
 section of \linebreak
$M \otimes (\otimes E^{*}) \otimes (\otimes E ).$\\

The notion of connections has been extended to Lie algebroids 
(see ~\cite{F} for example).
Let ${\cal B}$ be an ${\cal O}_{X}$-module. A $E$-connection on ${\cal B}$ is
a linear operator 
$$\nabla : \Gamma ({\cal B}) \to {}^{E}\Omega^{1}(\Gamma ({\cal B}))
=\Gamma \left ( Hom_{{\cal O}_{X}}(\wedge^{1}E, {\cal B})\right )$$
satisfying the following equation : for any $f \in \Gamma ({\cal O}_{X})$ and
any $v \in \Gamma ({\cal B})$ 
$$\nabla (fv)={}^{E}d(f)v + f \nabla (v).$$
If $u$ is an element of $E$, the connection $\nabla$ defines a map 
$\nabla_{u} : {\cal B} \to {\cal B}$. 

Assume now that ${\cal B}$ is a bundle. 
If $(e_{1}, \dots, e_{d})$ is a local basis of $E$ and 
$(b_{1}, \dots, b_{n})$ is a local basis of ${\cal B}$, one has 
$$\nabla_{e_{i}}(b_{j})=\Gamma_{i,j}^{k}b_{k}.$$  
The connection $\nabla$ is determined by its Christoffel symbol 
$\Gamma_{i,j}^{k}$. 

\begin{definition}
The curvature $R$ of a connection $\nabla$ with values in ${\cal B}$
is the section $R$ of the bundle 
$E^{*} \otimes E^{*}\otimes {\cal B}^{*}\otimes {\cal B}$ defined by :
for any $u,v$ in $\Gamma (E)$ and $b$ in $\Gamma ({\cal B})$ 
$$R(u,v)(b)=
\left ( \nabla_{u} \circ \nabla_{v} -\nabla_{v} \circ \nabla_{u}-
\nabla_{[u,v]}\right )(b)$$
\end{definition}
The curvature tensor is locally determined by the $(R_{i,j})_{k}^{l}$ 
 defined by 
$$R(e_{i},e_{j})b_{k}=(R_{i,j})_{k}^{l}b_{l}.$$
For a connection $\nabla$ on ${\cal B}=E$, one can define the torsion
tensor.
\begin{definition}
The torsion of $\nabla$ is a section of $E \otimes E^{*} \otimes E^{*}$ 
defined by~: for any $u$,$v$ in $\Gamma (E)$,
$$T(u,v)=\nabla_{u}(v)-\nabla_{v}(u)-[u,v].$$
\end{definition}
\begin{proposition}
A torsion free connection on $E$ exists.
\end{proposition}

A proof of this proposition can be found in ~\cite{C2}.\\

{\bf Examples of $D(E)$-modules}

{\it Example 1 :}

Flat connections provides examples of $D(E)$-modules.\\

{\it Example 2 :}
 
If $E$ is a Lie algebroid with anchor map $\omega$, then 
$Ker\omega $ is a left $D(E)$-module for the following operations :
for all $f$ in ${\cal O}_{X}$, for all $\xi$ in $E$, for all $\sigma$ in 
$Ker \omega$, 
$$f\cdot \sigma =f\sigma ,\;\;\; \xi \cdot \sigma =[\xi, \sigma ].$$

{\it Example 3 :}

If $M$ and $N$ are two left $D(E)$-modules, then 
(see ~\cite{Bo} for the
${\cal D}_{X}$-module case  and ~\cite{Ch2})
$M{\displaystyle \mathop \otimes_{{\cal O}_{X}}}N$ and 
${\cal H}om_{{\cal O}_{X}}(M,N)$, endowed with the
two operations described below, are left $D(E)$-modules : 
$$\begin{array}{l}
\forall m \in M, \forall n \in N,
\forall a \in {\cal O}_{X}, \forall \xi \in E , \\
a \cdot (m \otimes n) \cdot a = a \cdot m \otimes n\\
\xi \cdot (m \otimes n)  = \xi \cdot m \otimes n + m \otimes \xi \cdot n\\
\\
\forall \phi \in {\cal H}om_{{\cal O}_{X}}(M,N), \forall m \in M,
\forall a \in {\cal O}_{X}, \forall \xi \in E,\\
(a \cdot \phi )(m)=a\phi(m)\\ 
(\xi \cdot \phi ) (m)=\xi \cdot \phi(m)-\phi (\xi \cdot m).
\end{array}$$

{\it Example 4 :}

It is a well known fact (~\cite{Bo}, ~\cite{Ka}) that 
the ${\cal O}_{X}$-module of differential
forms of maximal degree, $\Omega^{dim X}_{X}$, is endowed with a right
${\cal D}_{X}$-module structure. We may extend this result 
 (~\cite{Ch1}) to $\Lambda^{d}(E^{*})$
where $d$ be the rank of $E$. Indeed $E$ acts on $\Lambda^{d}(E^{*})$ 
by the adjoint action. The action of an element $\xi$ of $E$ is called
the Lie derivative of $\xi$ and is denoted $L_{\xi}$. The 
${\cal  O}_{X}$-module  
$\Lambda^{d}(E^{*})$, endowed with the
following operations, 
$$\begin{array}{l}
\forall \sigma \in \Lambda^{d}(E^{*}), \forall \xi \in E, \forall f \in
{\cal O}_{X},\\
\sigma \cdot a = a\sigma \\
\sigma \cdot \xi=-L_{\xi}(\sigma)
\end{array}$$
is a right $D(E)$-module. \\

{\it Example 5 :}

If ${\cal M}$ and ${\cal N}$ are two right $D(E)$-modules, then 
${\cal H}om_{{\cal O}_{X}}({\cal M}, {\cal N})$, endowed with the two
following operations, 
$$\begin{array}{l}
\forall \phi \in {\cal H}om_{{\cal O}_{X}}({\cal M}, {\cal N}), \forall
m \in {\cal M}, \forall a \in {\cal O}_{X}, \forall \xi \in E , \\
(a \cdot \phi )(m)=\phi (m)\cdot a\\
(\xi \cdot \phi )(m)=-\phi (m)\cdot \xi + \phi (m \cdot \xi)
\end{array}$$
is a left $D(E)$-module (~\cite{Ch2}). 
This was already known for $D$-modules. 
In particular 
${\cal H}om_{{\cal O}_{X}}\left ( \Lambda^{d}(E^{*}), \Omega^{dim
    X}_{X}\right )$ is a left $D(E)$-module which is used in 
~\cite{ELW}
to define the modular class of $E$.\\

{\it Example 6 :}

If ${\cal M}$ is a right $D(E)$-module  and 
${\cal N}$ is a left  $D(E)$-module, then 
${\cal M}{\displaystyle \mathop \otimes_{{\cal O}_{X}}} {\cal N}$, 
endowed with the two following operations, 
$$\begin{array}{l}
\forall m \in {\cal M}, \forall n \in {\cal N},
\forall a \in {\cal O}_{X}, \forall \xi \in E , \\
(m \otimes n) \cdot a = m \otimes a \cdot n= m \cdot a \otimes n\\
(m \otimes n) \cdot \xi = m \cdot \xi \otimes n - m \otimes \xi \cdot n
\end{array}$$
is a right  $D(E)$-module (see ~\cite{Bo} for $D$-modules and 
~\cite{Ch2}). 
Given  any $D(E)$-module which is locally free of rank one, the
functor  
${\cal N}\mapsto {\cal E}{\displaystyle \mathop \otimes _{{\cal
      O}_{X}}}{\cal N}$ establishes an equivalence of categories
between left and right $D(E)$-modules. Its inverse functor is given by
${\cal M} \to {\cal H}om_{{\cal O}_{X}}({\cal E}, {\cal M})$. 
This equivalence of categories
is well known for $D$-modules (~\cite{Bo}, ~\cite{Ka})
and was generalized to Lie algebroids in ~\cite{Ch2}. 
In the case where $X=\R^{d}$ and $E=T\R^{d}$, 
this equivalence of categories is
particularly simple because we may choose 
$dx^{1}\wedge \dots \wedge dx^{d}$ as a basis of the 
${\cal O}_{R^{d}}$-module $\Omega_{X}^{d}$. There exists a unique
anti-isomorphism of 
${\cal D}_{\R^{d}}$, $\sigma$, such that $\sigma (f)=f$
and $\sigma ({\displaystyle \frac {\partial}{\partial x^{i}}})=
-{\displaystyle \frac {\partial}{\partial x^{i}}}$. Any left 
${\cal D}_{\R^{d}}$-module can be seen as a right 
${\cal  D}_{\R^{d}}$-module (and conversely) in the following way :
$$ \forall P \in {\cal D}_{\R^{d}}, \; \forall m \in M, \;\; 
m \cdot P= \sigma (P) \cdot m .$$

{\it Example 7 :}

 Let ${\cal D}b_{X}$ the sheaf of distributions over $X$. As 
${\cal O}_{X}$ is a left ${\cal D}_{X}$-module, 
${\cal D}b_{X}$ is a right ${\cal D}_{X}$-module (by transposition).   \\

{\it Example 8 :}

Let us recall our  definition of a Lie algebroid morphism 
(~\cite{Ch2})
which coincides with that of Almeida and Kumpera (~\cite{AK})

\begin{definition}
Let $\left ( E_{X}, \omega_{X} \right )$
 (respectively $\left ( E_{Y}, \omega_{Y} \right )$) be a Lie
algebroid over $X$ (respectively $Y$).
A morphism $\Phi$ from $\left ( E_{X}, \omega_{X} \right )$
to $\left ( E_{Y}, \omega_{Y} \right )$ is a pair $(f, F)$
such that
\begin{itemize}
\item
$f\; : \; X \to Y$ is a ${\cal C}^{\infty}$-morphism
\item
$F \; : \;E_{X} \to  f^{*} E_{Y}=
{\cal O}_{X}{\mathop \otimes_{f^{-1}{\cal O}_{Y}}} f^{-1}E_{Y}$
such that the two following conditions are satisfied:

1) the diagram

$$\xymatrix{
E_{X} \ar[d]_{\omega_{X}}^{}
\ar[r]^{F}
& f^{*}E_{Y}
\ar[d]_{}^{f^{*}\omega_{Y}}\\ 
\Theta_{X}  \ar[r]^{Tf}
&f^{*}\Theta_{Y}
}$$

commutes.

2)  Let $\xi$ and  $\eta$ be two elements of $E_{X}^{2}$.
Put  $F(\xi)={\displaystyle \sum _{i=1}^{m}}a_{i} \otimes \xi_{i}$ and
$F(\eta)={\displaystyle \sum _{j=1}^{m}}b_{j} \otimes \eta_{j}$, then
$$F([\xi, \eta]) =
{\displaystyle \sum _{j=1}^{n} } \omega_{X}(\xi)(b_{j}) \otimes \eta_{j}
-{\displaystyle \sum _{i=1}^{n} } \omega_{X}(\eta)(a_{i}) \otimes \xi_{i}
+ {\displaystyle \sum _{i,j}} a_{i}b_{j} \otimes [\xi_{i}, \eta_{j}]. $$

\end{itemize}
\end{definition}
If $\Phi=(f,F)$ is Lie algebroid morphism from $(E_{X}, \omega_{X})$ to 
$(E_{Y}, \omega_{Y})$ and ${\cal M}$ is a $D(E_{Y})$-module, then  
${\cal O}_{X}{\mathop \otimes _{f^{-1}{\cal O}_{Y}}}
f^{-1}{\cal M} $
 endowed with the two following operations

$$\begin{array}{l}
\forall (a,b) \in {\cal O}_{X}^{2}, \; \forall \xi \in E_{X},\;
\forall m \in f^{-1}{\cal M} \\
\\
a \cdot \left ( b  \otimes  m \right ) = a  b  \otimes  m \\
\xi \cdot \left ( b  \otimes  m \right ) =
\omega_{X}(\xi)(b) \otimes  m  +
\sum_{i} ba_{i} \otimes \xi_{i}m
\end{array}$$
(where $F(\xi)= {\displaystyle \sum_{i}}a_{i} \otimes \xi_{i}$
with $a_{i}$ in ${\cal O}_{X}$ and $\xi_{i}$ in $f^{-1}E_{Y}$) 
is a left $D(E_{X})$-module (~\cite{Ch2}). 

Morphisms of Lie algebroids generalize at the same time Lie algebra
morphisms and morphisms between ${\cal C}^{\infty}$-manifolds.  
Examples of Lie algebroid morphisms can be found in ~\cite{Ch3}. 
The $D(E_{X})\otimes f^{-1}D(E_{Y})^{op}$-module 
${\cal O}_{X}{\mathop \otimes _{f^{-1}{\cal O}_{Y}}}
f^{-1}D(E)$ generalizes the transfer module for $D$-modules 
(see ~\cite{Bo},~\cite{Ka}, ~\cite{Ch2})

\subsection{The sheaves of DGLAs $^{E}T_{poly}$ and $^{E}D_{poly}$}

The sheaf of DGLAs of polyvectorfields can be  extended to the Lie algebroids
setting. 
The sheaf of DGLAs  $^{E}T_{poly}$ of $E$-polyvector fields 
is defined as follows (~\cite{C1}):
$$^{E}T_{poly}=
{\displaystyle \mathop \oplus_{k \geq -1}}  ^{E}T^{k}_{poly}
= {\displaystyle \mathop \oplus_{k \geq -1}  \wedge^{k+1}E }$$
endowed with the zero differential and the Lie bracket 
$[\;, \; ]_{S}$ uniquely
defined by the two following properties : 
$$\begin{array}{l}
\bullet \forall f,g \in {\cal O}_{X}, [f,g]_{S}=0\\
\bullet \forall \xi \in E, \; \forall f \in {\cal O}_{X}, \; 
[\xi, f ]_{S}=\omega (\xi)(f)\\
\bullet \forall \xi, \eta \in E, \;\; [\xi, \eta ]_{S}=[\xi, \eta]_{E}\\
\bullet \forall u \in {}^{E}T_{poly}^{k}, v \in {}^{E}T_{poly}^{l}, w \in
{}^{E}T_{poly},\\
\; [u, v \wedge w ]_{S}=[ u,v]_{S}\wedge w + (-1)^{k(l+1)}v \wedge [ u,w]_{S}
\end{array}$$
In ~\cite{C1}, 
D. Calaque extended the sheaf of DGLAs of polydifferential operators
to the Lie algebroid setting.  Before recalling his construction, let
us fix some notations. \\

{\bf Notation :}

Let $M_{0}, M_{1}, \dots , M_{n}$ be $n+1$ $D(E)$-modules. Denote by 
$\pi_{i}: D(E) \to End(M_{i})$ the maps defined by these actions. 
An element $P_{0} 
{\displaystyle \mathop \otimes _{{\cal O}_{X}}}\dots 
{\displaystyle \mathop \otimes _{{\cal O}_{X}}}P_{n}$ of 
$D(E)^{\otimes n+1}$ defines a map 
$$\begin{array}{rcl}
\pi _{0}(P_{0}) \otimes \dots \otimes \pi _{n+1}(P_{n+1}) : 
M_{0} {\displaystyle \mathop \otimes _{\R_{X}}}
\dots {\displaystyle \mathop \otimes _{\R_{X}}}M_{n}& \to&
M_{0} {\displaystyle \mathop \otimes _{{\cal O}_{X}}}
\dots {\displaystyle \mathop \otimes _{{\cal O}_{X}}}M_{n}\\
m_{0}  {\displaystyle \mathop \otimes _{\R}}
\dots {\displaystyle \mathop \otimes _{\R}}m_{n}& \mapsto &
\pi_{0}(P_{0})(m_{0})
{\displaystyle \mathop \otimes _{{\cal O}_{X}}}\dots 
{\displaystyle \mathop \otimes _{{\cal O}_{X}}}
\pi_{n}(P_{n})(m_{n}).\\
\end{array}$$
In the sequel, we will be in the following situation 
$M_{0}, \dots, M_{i-1}, M_{i+1}, \dots, M_{n}$ are $D(E)$ endowed with
left multiplication. If $P$ is in $D(E)$, we will then write $P$ for
left multiplication with $P$ which amounts to omit 
$\pi_{i}$. The $D(E)$-module $M_{i}$ will be ${\cal O}_{X}$ (with its
natural $D(E)$-module structure) and we
will write $\omega$ (as the anchor map) for the map from $D(E)$ to 
$End({\cal O}_{X})$.  \\

Calaque defines the sheaf of DGLAs $^{E}D_{poly}^{*}$ of $E$-polydifferential
operators as follows : 
$$^{E}D_{poly}^{*}=
{\displaystyle \mathop \oplus_{k \geq -1}}  ^{E}D^{k}_{poly}$$ 
where 
$$\begin{array}{l}
^{E}D^{-1}_{poly}={\cal O}_{X} \\
^{E}D^{k}_{poly}= D(E)^{{\displaystyle \mathop \otimes}^{k+1}_{{\cal
      O}_{X}}}
\;\;{\rm if}\;\;k \geq 0.\\
\end{array}$$
Before defining the Lie bracket over $^{E}D_{poly}^{*}$, we need to
introduce the bilinear product of degree $0$, 
$$\bullet : {}^{E}D_{poly}^{*} \otimes {}^{E}D_{poly}^{*} 
\to {}^{E}D_{poly}^{*}.$$
Let $P$ (respectively $Q$) be an homogeneous element of $^{E}D_{poly}^{*}$
of positive degree $\mid P \mid$ (respectively $\mid Q \mid $), and 
let $f$ (respectively $g$) be an element of 
$^{E}D_{poly}^{-1}={\cal O}_{X}$. We have :
$$\begin{array}{l}
P \bullet Q = {\displaystyle \mathop \sum _{i=0}^{\mid P \mid}}(-1)^{i \mid Q \mid}
\left (id^{\otimes^{i}}\otimes \Delta^{(\mid Q \mid)}\otimes 
id^{\otimes \mid P \mid  -i}\right )(P)\cdot
(1^{\otimes^{i}}
{\displaystyle \mathop \otimes_{\R}} Q
 {\displaystyle \mathop \otimes_{\R}} 
1^{\otimes \mid P \mid -i})\\
P \bullet f = {\displaystyle \mathop \sum _{i=0}^{\mid P \mid}}(-1)^{i}
\left ( id^{\otimes^{i}}\otimes \omega \otimes id^{\otimes \mid P \mid
  -i}\right )(P)\cdot
(1^{\otimes^{i}}
{\displaystyle \mathop \otimes_{\R}} f 
{\displaystyle \mathop \otimes_{\R}} 
1^{\otimes \mid P \mid -i})\\
f\bullet g=0\\
f \bullet P =0
\end{array}$$
The Lie bracket between $P_{1}\in {}^{E}D_{poly}^{k_{1}}$ and 
$P_{2} \in {}^{E}D_{poly}^{k_{2}}$ is 
$$[P_{1},P_{2}]=P_{1}\bullet P_{2}-(-1)^{k_{1}k_{2}}P_{2}\bullet P_{1}.$$
The differential on $^{E}D_{poly}$ is $\partial =[1 \otimes 1, -]$.\\

Calaque has proved the following theorem (~\cite{C1}) which
generalizes Kontsevitch's result (~\cite{Ko})

\begin{theorem}
There exists a quasi-isomorphism of 
$L_{\infty}$-algebras, ${\Upsilon}$, from
$\Gamma \left (^{E}T_{poly}^{*}\right )$ to 
$\Gamma \left ( ^{E}D^{*}_{poly}\right )$. 
In other words, $\Gamma \left ( {}^{E}D_{poly}^{*}\right )$
is formal. 
\end{theorem}

\section{Main results}

Let $E$ be a Lie algebroid over a manifold $X$ and let $D(E)$ be the
sheaf of $E$- differential operators. We will
denote by $M$  a left $D(E)$-module.  

\subsection{The $^{E}T_{poly}^{*}$-module $^{E}T_{poly}^{*}(M)$.}

We introduce the complex $^{E}T_{poly}^{*}(M)$ of $E$-polyvector fields
with values in $M$
$$^{E}T_{poly}^{*}(M)=
{\displaystyle \mathop \oplus_{k \geq -1} {}^{E}T^{k}_{poly}(M)}
= {\displaystyle \mathop \oplus_{k \geq -1}  
\wedge^{k+1}E \otimes M }$$
 with differential zero. If $m$ is in $M$, we will identify $m$ with 
$1 \otimes m$. 

\begin{proposition}
$^{E}T_{poly}^{*}(M)$ is endowed with a $^{E}T^{*}_{poly}$-module structure
described as follows : 
for all $u=\xi_{1}\wedge \dots \wedge \xi_{k+1} \in
{}^{E}T^{k}_{poly}$ , 
$v \in  ^{E}T^{l}_{poly}$ (with $k,l \geq 0$), 
$f \in {\cal O}_{X}$, $m\in M$, 
$$\begin{array}{l} 
\bullet f \cdot _{S} m = 0\\
\bullet 
(\xi_{1}\wedge \dots \wedge \xi_{k+1}) \cdot _{S}m =
{\displaystyle \mathop \sum _{i=1}^{k+1}}(-1)^{k+1-i}
 \xi_{1}\wedge \dots \wedge \widehat{\xi_{i}}\wedge \dots \wedge \xi_{k+1}
\otimes \xi_{i}\cdot m \\
\bullet f \cdot_{S} (v \otimes m)=[f,v]_{S} \otimes m \\
\bullet 
u \cdot_{S} (v \otimes m )= 
[u,v]_{S} \otimes m + (-1)^{k(l+1)}v \wedge u\cdot_{S} m .\\ 
\end{array}$$
\end{proposition}

When there is no ambiguity, we will drop the subscript $S$ in the
notation of the action of $^{E}T_{poly}^{*}$ over $^{E}T_{poly}^{*}(M)$.\\

{\it Proof of the proposition}

It is easy to check that the actions above are well defined. 
Let  $a$ be in  $^{E}T^{s}_{poly}$.
We need to verify that the following relation holds 
$$ u \cdot \left ( v \cdot \left ( a \otimes m \right )\right )
-(-1)^{kl}
v \cdot \left ( u \cdot \left ( a \otimes m \right )\right)=
[u,v]\cdot \left (a \otimes m \right ).$$
A straightforward computation shows that it is enough to check this
relation for $a=1$, which we will assume.We will need the two
following lemmas.
\begin{lemma}
if $a \in {}^{E}T^{*}_{poly}$, $u \in {}^{E}T^{k}_{poly}$, 
$v \in {}^{E}T^{l}_{poly}$ ($k,l \geq -1$),  one has 
$$u \cdot \left ( v \wedge a \otimes m \right ) = 
[u,v] \wedge a \otimes m + 
(-1)^{k(l+1)} v \wedge u \cdot \left ( a \otimes m \right ).$$ 
\end{lemma}
 
{\it Proof of the lemma :} It is a straightforward
  computation. $\Box$

\begin{lemma}
Let $a \in {}^{E}T^{*}_{poly}$, $m \in M$, $k,l \geq 0$, 
$u \in {}^{E}T^{k}_{poly}$, 
$v\in {}^{E}T^{l}_{poly}$. One has the following relation 
$$(u \wedge v) \cdot (a \otimes m) = 
u \wedge \left ( v \cdot (a \otimes m) \right ) + 
(-1)^{(k+1)(l+1)} v \wedge \left ( u \cdot (a \otimes m)\right )$$
\end{lemma} 

{\it Proof of the lemma :}

An easy computation shows that we may assume $a=1$.
The proof of the lemma goes by induction over $k$.
The case $k=0$ is obvious so that we assume $k \geq 1$. 
Set $u=\xi_{1} \wedge \dots \wedge \xi_{k+1}$ and 
$u'=\xi_{2} \wedge \dots \wedge \xi_{k+1}$ 
so that $u=\xi_{1}\wedge u'$. 
Using the induction hypothesis and the case $k=0$, we get the
following sequence of equalities.
$$\begin{array}{l}
(u\wedge v )\cdot m\\ 
=(-1)^{l+k+1}(u' \wedge v) \otimes \xi_{1}\cdot m 
+ \xi_{1} \wedge \left ( (u' \wedge v ) \cdot m \right ) \\
=(-1)^{l+k+1+k(l+1)} v \wedge u' \otimes \xi_{1}\cdot m
+ \xi_{1} \wedge u' \wedge (v \cdot m) 
+(-1)^{k(l+1)}\xi_{1}\wedge v \wedge (u' \cdot m) \\
= u \wedge (v \cdot m) + (-1)^{(k+1)(l+1)}v \wedge (u \cdot m). \Box \\ 
\end{array}$$

We will show the relation 
$$ u \cdot \left ( v \cdot  m \right )
-(-1)^{kl}
v \cdot \left ( u \cdot m \right)=
[u,v]\cdot  m $$
by induction on $l$.\\

{\it First case : l=-1}\\

In this case $v$ is a function on $X$ which will be denoted $f$. 
We proceed by induction over $k$. The cases $k=-1$ or $k=0$ are
obvious so that we assume $k\geq 1$. 
We set $u=\xi_{1} \wedge \dots \wedge \xi_{k+1}$ and 
$u'=\xi_{2} \wedge \dots \wedge \xi_{k+1}$. 

Using the two previous lemmas and the induction hypothesis , 
we get the following sequence of
equalities  :
$$\begin{array}{l}
u\cdot (f \cdot m)- (-1)^{k}f \cdot \left ( u \cdot m \right ) \\
= - (-1)^{k}f \cdot \left ( u \cdot m \right ) \\
= -(-1)^{k}f \cdot \left ( 
\xi_{1}\wedge \left ( u'\cdot m \right ) +
(-1)^{k}u' \otimes \xi_{1} \cdot m \right )\\
 = - (-1)^{k} [f, \xi_{1}] \left (u'\cdot m\right ) + 
(-1)^{k} \xi_{1} \wedge 
\left (f \cdot \left ( u'\cdot m \right )\right ) 
- [f, u'] \otimes \xi_{1}\cdot m\\
 =- (-1)^{k} [f, \xi_{1}]
 \left (u'\cdot m\right ) + 
(-1)^{k} \xi_{1} \wedge \left ( [f, u']\cdot m \right ) 
- [f, u'] \otimes \xi_{1}\cdot m\\
\end{array}$$
On the hand, 
$$[f, u]= [f, \xi_{1}] u' - \xi_{1}\wedge [f, u'].$$
hence, 
$$\begin{array}{l}
[f, u]\cdot m = [f, \xi_{1}] u'\cdot m 
-\xi_{1}\wedge \left ( [f,u'] \cdot m \right )
-(-1)^{k+1}[f,u']\otimes \xi_{1}\cdot m. 
\end{array}$$
The cas $l=-1$ follows.\\

{\it Second case : l=0}\\

In this case $v$ is an element of $E$ which will be denoted $\eta$. 
We proceed by induction over $k$. The cases $k=-1$ or $k=0$ are
obvious so that we assume $k\geq 1$. 
We set $u=\xi_{1} \wedge \dots \wedge \xi_{k+1}$ and 
$u'=\xi_{2} \wedge \dots \wedge \xi_{k+1}$. 

Using the two previous lemmas, we get the following sequence of
equalities  :
$$\begin{array}{l}
u\cdot (\eta \cdot m)- \eta \cdot \left ( u \cdot m \right ) \\
=\xi_{1} \wedge \left ( u' \cdot (\eta \cdot m)\right )
+(-1)^{k}u'\otimes \xi_{1}\cdot (\eta \cdot m )
-\eta \cdot \left (\xi_{1} \wedge (u' \cdot m) 
+(-1)^{k}u'\otimes \xi_{1}\cdot m \right )\\
=\xi_{1}\wedge \left ( [u', \eta] \cdot m \right ) + 
(-1)^{k} u' \otimes [\xi_{1},\eta]\cdot m 
-[\eta, \xi_{1}] \wedge (u' \cdot m) 
-(-1)^{k}[\eta , u']\otimes \xi_{1}\cdot m .
\end{array}$$
On the other hand, 
$$[u, \eta ]= -[\eta ,\xi_{1}] \wedge u' - \xi_{1}\wedge [\eta, u'].$$
hence
$$[u, \eta]\cdot m=
-[\eta,\xi_{1}]\wedge (u' \cdot m)-
(-1)^{k}u' \otimes [\eta, \xi_{1}]\cdot m
-(-1)^{k}[\eta , u']\otimes \xi_{1} \cdot m 
-\xi_{1}\wedge [\eta, u']\cdot m.$$

{\it Third case : $l\geq 1$}\\

We proceed by induction.  
We set $v=\eta_{1} \wedge \dots \wedge \eta_{k+1}$ and 
$u'=\eta_{2} \wedge \dots \wedge \eta_{k+1}$. 
Using the previous lemmas and the induction hypothesis, 
we get the following sequences of
equalities :
$$\begin{array}{l}
u \cdot \left ( v \cdot  m \right )-(-1)^{kl} 
v \cdot \left ( u \cdot m \right )\\
= u \cdot \left ( \eta_{1} \wedge (v' \cdot m) +
(-1)^{l} v' \otimes \eta_{1}\cdot m \right )
-(-1)^{kl}\eta_{1}\wedge \left ( v' \cdot (u \cdot m)\right )\\
-(-1)^{lk+l}v' \wedge \left ( \eta_{1} \cdot (u \cdot m) \right )\\
=(-1)^{k(l+1)}v'\wedge \left ( [u,\eta_{1}]\cdot m \right )
+(-1)^{k}\eta_{1}\wedge [u,v']\cdot m 
+[u, \eta_{1}] \wedge (v'\cdot m)\\
+(-1)^{l} [u,v'] \otimes \eta_{1}\cdot m .\\
\end{array}$$
On the other hand,

$$[u,v]=[u, \eta_{1}]\wedge v'  + (-1)^{k}\eta_{1}\wedge [u,v']$$
hence
$$\begin{array}{rcl}
[u,v]\cdot m &=&
[u, \eta_{1}]\wedge (v'\cdot m) +
(-1)^{l(k+1)} v'\wedge \left ( [u, \eta_{1}]\cdot m \right )
+(-1)^{k}\eta_{1}\wedge \left ( [u,v']\cdot m \right )\\
&+&(-1)^{l}[u,v']\otimes \eta_{1}\cdot m.\\
\end{array}$$
The case $l\geq 1$ follows. $\Box$

\subsection{The $^{E}D_{poly}^{*}$-module $^{E}D_{poly}^{*}(M)$}

Let $M$ be a $D(E)$-module. 
Denote by $\pi$ the map from $D(E)$ to $End(M)$ \linebreak 
determined by the left  
$D(E)$-module structure on $M$. We will use the same 
\linebreak notation as in
section 2.2. 
We will also use the map $\tau_{i}$ from \linebreak  
$\left ( {\displaystyle \mathop \otimes^{i-1}_{{\cal O}_{X}}}D(E)
\right )
{\displaystyle \mathop \otimes_{{\cal O}_{X}}} 
\left ( D(E){\displaystyle \mathop \otimes_{{\cal O}_{X}}M} \right )
{\displaystyle \mathop \otimes_{{\cal O}_{X}}}
\left ( {\displaystyle \mathop \otimes^{q+1-i}_{{\cal O}_{X}}}D(E)
\right )$ to 
$D(E)^{{\displaystyle \mathop \otimes}^{q+1}_{{\cal O}_{X}}}
{\displaystyle \mathop \otimes_{{\cal O}_{X}}} M$
defined by 
$$\tau_{i} \left ( Q_{0}\otimes \dots \otimes Q_{i-1} \otimes
  (Q_{i}\otimes m) \otimes Q_{i+1}\otimes \dots \otimes Q_{q} \right )
=Q_{0}\otimes \dots \otimes Q_{q} \otimes m. $$
Let us introduce the complex  $^{E}D_{poly}(M)$ of $E$-polydifferential
operators with values in $M$ as follows : 
$$^{E}D_{poly}(M)=
{\displaystyle \mathop \oplus_{k \geq -1}}  ^{E}D^{k}_{poly}(M)$$ 
where 
$$\begin{array}{l}
^{E}D^{-1}_{poly}(M)= M \\
^{E}D^{k}_{poly}(M)= 
D(E)^{{\displaystyle \mathop \otimes}^{k+1}_{{\cal O}_{X}}}
{\displaystyle \mathop \otimes_{{\cal O}_{X}}} M
\;\;{\rm if}\;\;k \geq 0.\\
\end{array}$$ 
Let us define two maps denoted in the same way 
$$\begin{array}{l}
\bullet : {}^{E}D_{poly}^{*} \otimes {}^{E}D_{poly}^{*}(M) 
\to {}^{E}D_{poly}^{*}(M)\\
\bullet : {} ^{E}D_{poly}^{*}(M)\otimes {}^{E}D_{poly}^{*} 
\to {}^{E}D_{poly}^{*}(M).
\end{array}$$
If $P$ and $Q$ are homogeneous elements of $^{E}D_{poly}^{*}$ of
non negative degree respectively $\mid P \mid$ and $\mid Q \mid$, 
$f$ is an element of $^{E}D_{poly}^{-1}$ 
and  $m$ is in $M$, then 
$$\begin{array}{l}
P \bullet (Q \otimes m)= {\displaystyle \mathop \sum _{i=0}^{\mid P \mid}}
(-1)^{i \mid Q \mid} \tau_{i} \left [
\left (id^{\otimes^{i}}\otimes \Delta^{(\mid Q \mid+1)}\otimes 
id^{\otimes \mid P \mid  -i}\right )(P)\cdot
(1^{\otimes^{i}}
{\displaystyle \mathop \otimes_{\R}} 
(Q \otimes m) 
{\displaystyle \mathop \otimes_{\R}} 
1^{\otimes \mid P \mid -i})
\right ]\\
P \bullet m = \tau_{i} \left [
{\displaystyle \mathop \sum _{i=0}^{\mid P \mid}}(-1)^{i}
\left ( id^{\otimes^{i}}\otimes \pi \otimes id^{\otimes \mid P \mid
  -i}\right )(P)\cdot
(1^{\otimes^{i}}
{\displaystyle \mathop \otimes_{\R}} m 
{\displaystyle \mathop \otimes_{\R}} 
1^{\otimes \mid P \mid -i}) \right ] \\
f\bullet m=0\\
f \bullet (Q \otimes m) =0\\
(Q \otimes m) \bullet P = Q \bullet P \otimes m \\
m \bullet P =0\\
m \bullet f =0.\\
\end{array}$$
Note that the second, the third and the fourth equations could be recovered
 from the first one. 
The differential, $\partial_{M}$, on $^{E}D_{poly}^{*}(M)$ is given by :
for all $Q \otimes m$ in $^{E}D_{poly}^{*}(M)$,
$$\begin{array}{rcl} 
\partial _{M}(Q \otimes m) & = & (1\otimes 1) \bullet (Q \otimes m) -
(-1)^{\mid Q \mid}(Q \otimes m) \bullet (1\otimes1)\\
&=& \partial (Q)\otimes m\\
\end{array}$$
where $1 \otimes 1 \in {}^{E}D^{1}_{poly}$.
\begin{theorem}
$^{E}D_{poly}^{*}(M)$ is  endowed with a $^{E}D^{*}_{poly}$-module
structure as follows. 
$$\begin{array}{l}
\forall P \in {}^{E}D^{p}_{poly}, 
\forall (Q \otimes m) \in {}^{E}D^{q}_{poly}(M)\\
P \cdot_{G} (Q \otimes m)= P \bullet (Q \otimes m) -
(-1)^{pq} (Q\otimes m) \bullet P
\end{array}$$
\end{theorem}
{\it Proof of the theorem :}

Let $P \in {}^{E}D_{poly}^{p}$,  $Q \in {}^{E}D_{poly}^{q}$, 
 $\lambda \in {}^{E}D_{poly}^{r}(M)$. 
Introduce the following quantity 
$$A(P,Q, \lambda )=(P \bullet Q )\bullet \lambda - 
P \bullet ( Q \bullet \lambda).$$
The theorem follows from the  lemma below.
\begin{lemma}
The following equality holds : 
$$A(P, Q, \lambda)=(-1)^{qr}A(P, \lambda, Q)$$
\end{lemma}

This lemma is well known in the case where $E=TX$ and 
$M={\cal O}_{X}$ 
(see for example the article of Keller in ~\cite{BCKT}).

In the general case, it  follows from a  straightforward but 
tedious computation. $\Box$

\subsection{The Hochschid-Kostant-Rosenberg theorem}

\begin{theorem}
The map $U_{HKR}^{M}$ from $\left ( ^{E}T_{poly}^{*}(M), 0 \right )$ to 
$\left ( ^{E}D_{poly}^{*}(M), \partial_{M} \right )$ defined \linebreak 
by : for all 
$v_{1}, \dots , v_{n}$ in $E$ and all $m$ in $M$,  
$$\begin{array}{l}
U_{HKR}^{M}(v_{0} \wedge \dots \wedge v_{n} \otimes m)=
{\displaystyle \frac{1}{(n+1)!}}
{\displaystyle \sum_{\sigma \in S_{n+1}}}\epsilon (\sigma)
v_{\sigma (0)}\otimes \dots \otimes v_{\sigma(n)}\otimes m \\
U_{HKR}^{M}(m)=m
\end{array}$$
is a quasi-isomorphism.
\end{theorem}

The first one to have proved such a statement in the affine case 
(for $E=TX$ and $M={\cal O}_{X}$) seems to be J. Vey (~\cite{V}). 
A proof for the tangent bundle
of any manifold (and $M={\cal O}_{X}$) can be found in ~\cite{Ko}. 
This theorem is proved in ~\cite{C1} for any Lie algebroid and 
$M={{\cal O}_{X}}$. \\

{\it Proof of the theorem :}

This theorem will be  a consequence of the proof of theorem 
~\ref{formality} and of the 
following well known result.

\begin{lemma}
  
$T$ be a finite dimensional $\R$-vector space. Consider
the complex 
$\wedge^{*}T
={\displaystyle \mathop \oplus _{p \in \N}}\wedge ^{p} T$
 with zero differential and the complex 
${\displaystyle \mathop \oplus _{p \in \N}}
\left ( {\displaystyle \mathop \otimes ^{p}}S(E)\right )$
with the differential  
$$\partial =  
id^{\otimes p}\otimes 1 + (-1)^{p-1}1 \otimes id^{\otimes p}
+ (-1)^{p-1}
{\displaystyle \mathop \sum_{i=0}^{n}}(-1)^{i}
 id^{\otimes i}\otimes \Delta \otimes id^{\otimes n-i} .$$ 
The $\R$-linear map ${\Theta}$ from 
$\bigwedge^{*} T $ to 
${\displaystyle \mathop \oplus _{p \in \N}}
 {\displaystyle \mathop \otimes ^{p}}S(T)$ defined by : for all 
$v_{1}, \dots , v_{p}$ in $T$,  
$$\begin{array}{l}
{\Theta}(v_{0} \wedge \dots \wedge v_{p})=
{\displaystyle \frac{1}{(p+1)!}}
{\displaystyle \sum_{\sigma \in S_{p+1}}}\epsilon (\sigma)
v_{\sigma (0)}\otimes \dots \otimes v_{\sigma(p)} \\
{\Theta}(1)=1
\end{array}$$
is a quasi-isomorphism.
\end{lemma}

\subsection{Main statement}

We have seen that $\Gamma \left (^{E}D_{poly}^{*}(M)\right )$ 
is a module over the DGLA $\Gamma \left (^{E}D_{poly}^{*}\right )$. 
As we know (~\cite{C1}) that there is a $L_{\infty}$-morphism
from $\Gamma \left (^{E}T_{poly}^{*}\right )$ to 
$\Gamma \left ( ^{E}D_{poly}^{*}\right )$, 
we deduce that
$\Gamma \left ( ^{E}D_{poly}^{*}(M)\right )$ 
is naturally endowed with the structure of an 
$L_{\infty}$-module over the DGLA 
$\Gamma \left (^{E}T_{poly}^{*}\right )$. 
We can now state the main result of this paper.
\begin{theorem}\label{formality}
There is a quasi-isomorphism of $L_{\infty}$-modules 
over  $\Gamma \left (^{E}T_{poly}^{*}\right )$ from 
$\Gamma \left (^{E}T_{poly}^{*}(M)\right )$ to  
$\Gamma \left (^{E}D_{poly}^{*}(M)\right )$ 
that induces $U^{M}_{HKR}$ in cohomology.
\end{theorem}
Our result extends Calaque 's result 
(~\cite{C1}, take $M={\cal O}_{X}$) and 
Kontsevitch 's result (~\cite{Ko}, take $M={\cal O}_{X}$ and $E=TX$).

\section{Proof}

The proof is analogous to that of 
~\cite{D1}, ~\cite{C1}, ~\cite{D2}, ~\cite{CDH}. 

\subsection{Fedosov Resolutions}
As before, $E$ will denote a Lie algebroid and $M$ will be a 
$D(E)$-module. 

Following Fedosov and Dolgushev (~\cite{Fe}, ~\cite{D1}), 
Calaque introduced (~\cite{C1}, see also ~\cite{CDH}), 
the locally free ${\cal O}_{X}$-modules ${\cal W}=\widehat{S}(E^{*})$, 
${\cal  T}^{*}$ and ${\cal D}^{*}$. Let
us recall their definition. 

${\bullet}$ ${\cal W}=\widehat{S}(E^{*})$ is the 
locally free ${\cal O}_{X}$-module  whose sections
are functions that are formal in the fiber. An element $s$ of 
$\Gamma (U, {\cal W})$ can be locally written 
$$s={\displaystyle \mathop \sum_{l=0}^{\infty}}
s_{i_{1}, \dots, i_{l}}y^{i_{1}}\dots y^{i_{l}}$$
where $y^{1}, \dots, y^{d}$ are coordinates in the fiber of $E$ and 
$s_{i_{1}, \dots , i_{l}}$ are coefficients of a symmetric
covariant $E$-tensor.\\

${\bullet}$ 
${\cal T}^{*}={\cal W}{\displaystyle \mathop \otimes_{{\cal O}_{X}}} 
\bigwedge ^{*+1}E$ is
the graded  locally free ${\cal O}_{X}$-module 
of formal fiberwise polyvector fields on $E$ with
shifted degree. 
A homogeneous section of degree $k$ of ${\cal T}^{*}$ can
be locally written 
$${\displaystyle \mathop \sum_{l=0}^{\infty}
v_{i_{1},\dots, i_{l}}^{j_{0},\dots,j_{k}}}
y^{i_{1}} \dots y^{i_{l}}
{\displaystyle \frac{\partial}{\partial y^{j_{0}}}}
\wedge \dots \wedge  
{\displaystyle \frac{\partial}{\partial y^{j_{k}}}}$$
where $v_{i_{1},\dots, i_{l}}^{j_{0},\dots,j_{k}}$ 
are components of an $E$-tensor 
symmetric covariant in the indices $i_{1}, \dots,i_{l}$, contravariant 
antisymmetric in the indices $j_{0}, \dots, j_{k}$ 

$\bullet$ ${\cal D}^{*}=\widehat{S}(E^{*})
{\displaystyle \mathop \otimes_{{\cal O}_{X}}} 
T^{*+1}\left (S(E)\right )$ is the graded
 locally free ${\cal O}_{X}$-module of formal fiberwise
$E$-polydifferential operators with shifted degree. An homogeneous
section of degree $k$ of ${\cal D}^{*}$ can be locally written 
$${\displaystyle \mathop \sum_{l=0}^{\infty}}
P_{i_{1}, \dots, i_{l}}^{\alpha_{0}, \dots, \alpha_{k}}(x)
y^{i_{1}}\dots y^{i_{l}}
{\displaystyle \frac{\partial^{\mid \alpha_{0}\mid}}
{\partial y^{\alpha_{0}}}} \otimes \dots \otimes 
{\displaystyle \frac{\partial^{\mid \alpha_{k}\mid}}
{\partial y^{\alpha_{k}}}}$$
where the $\alpha_{i}$'s are multi-indices, 
the $P_{i_{1}, \dots, i_{l}}^{\alpha_{0}, \dots, \alpha_{k}}(x)$
 are components of an $E$-tensor with obvious symmetry.\\

We will need to introduce the ${\cal O}_{X}$-modules 
${\cal D}^{*}(M)$ and ${\cal T}^{*}(M)$. 

${\bullet}$ ${\cal T}^{*}(M)$ is
the graded ${\cal O}_{X}$-module  of formal fiberwise polyvector
fields on $E$ with
values in $M$ with shifted degree. 
A homogeneous section of degree $k$ of ${\cal T}^{*}(M)$ can
be locally written 
$${\displaystyle \mathop \sum_{l=0}^{\infty}
m_{i_{1},\dots, i_{l}}^{j_{0},\dots,j_{k}}}
y^{i_{1}} \dots y^{i_{l}}
{\displaystyle \frac{\partial}{\partial y^{j_{0}}}}
\wedge \dots \wedge  
{\displaystyle \frac{\partial}{\partial y^{j_{k}}}}$$
where $m_{i_{1},\dots, i_{l}}^{j_{0},\dots,j_{k}}$ 
are components of an $E$-tensor with values in $M$
symmetric covariant in the indices $i_{1}, \dots,i_{l}$, contravariant 
antisymmetric in the indices $j_{0}, \dots, j_{k}.$\\

${\bullet}$
 ${\cal D}^{*}(M)$ 
is the graded ${\cal O}_{X}$-modules of formal fiberwise
$E$-polydifferential operators with values in $M$ 
(with shifted degree). A homogeneous
section of degree $k$ of ${\cal D}^{*}(M)$ can be locally written 
$${\displaystyle \mathop \sum_{l=0}^{\infty}}
\mu_{i_{1}, \dots, i_{l}}^{\alpha_{0}, \dots, \alpha_{k}}(x)
y^{i_{1}}\dots y^{i_{l}}
{\displaystyle \frac{\partial^{\mid \alpha_{0}\mid}}
{\partial y^{\alpha_{0}}}} \otimes \dots \otimes 
{\displaystyle \frac{\partial^{\mid \alpha_{k}\mid}}
{\partial y^{\alpha_{k}}}}$$
where the $\alpha_{i}$'s are multi-indices, 
the $\mu_{i_{1}, \dots, i_{l}}^{\alpha_{0}, \dots, \alpha_{k}}(x)$
 are coefficients of an $E$-tensor with values in $M$
with obvious symmetry.\\

{\bf Remark :}

One has the obvious equality 
${\cal T}^{*}({\cal O}_{X})={\cal T}^{*}$ and 
${\cal D}^{*}({\cal O}_{X})={\cal D}^{*}$. \\

\noindent {\bf Notation :}
 
Let $\R ^{d}_{formal}$ be the formal completion of 
$\R ^{d}$ at the
origin. The ring of functions on $\R^{d}_{formal}$ is 
$\R [[y^{1}, \dots, y^{d}]]$ and the Lie-Rinehart algebra of vector
fields is $Der \left ( \R [[y^{1}, \dots, y^{d}]] \right )$. 
Denote by 
$T_{poly}^{*}(\R ^{d}_{formal})$ and 
$D_{poly}^{*}(\R^{d}_{formal})$ the 
DGLAs of polyvector fields and polydifferential
operators on  $\R ^{d}_{formal}$ respectively. If 
$t_{1} \in D_{poly}^{k_{1}-1}(\R^{d}_{formal})$ and 
$t_{2} \in D_{poly}^{k_{2}-1}(\R^{d}_{formal})$, one defines
their cup-product 
$t_{1}\sqcup t_{2}\in
D_{poly}^{k_{1}+k_{2}-1}(\R^{d}_{formal})$ by :
$$\begin{array}{l}
\forall a_{1}, \dots, a_{k_{1}+k_{2}} \in \R [[y^{1}, \dots, y^{d}]],\\
(t_{1}\sqcup t_{2})(a_{1}, \dots, a_{k_{1}+k_{2}})=
t_{1}(a_{1}, \dots, a_{k_{1}}) t_{2}(a_{k_{1}+1}, \dots, 
a_{k_{1}+k_{2}})\\
\end{array}$$
The cup-product endows  $D_{poly}^{*}(\R^{d}_{formal})$ with
the structure of a DGAA.\\

{\bf Remark: }

Fiberwise product endows ${\cal W}$ with the structure of 
bundle of commutative algebra.  ${\cal T}^{*}$ is a
differential Lie algebra with zero differential and Lie bracket
induced by fiberwise Schouten bracket on 
$T_{poly}(\R ^{d}_{formal})$. Similarly, fiberwise Schouten bracket
allows to endow   ${\cal T}^{*}(M)$ with a 
${\cal  T}^{*}$-module structure. We can make the same type of remark
for ${\cal D}$, ${\cal D}(M)$ and the Gerstenhaber bracket.\\

Let ${\cal B}$ be any  of  the ${\cal O}_{X}$-modules  
introduced above. We will need
to tensor ${\cal B}$ by $\bigwedge ^{*} (E^{*})$. We set 
$^{E}\Omega ({\cal B})= \bigwedge ^{*} (E^{*}) \otimes {\cal B}$.\\

{\bf Structures on $^{E}\Omega \left ({\cal B} \right )$}

${\bullet}$
$^{E}\Omega ({\cal W})$ is a bundle of graded commutative algebras 
with grading given by exterior degree of $E$-forms.

${\bullet}$
Schouten bracket on $T^{*}_{poly}(\R^{d}_{formal})$ induces a structure of
sheaf of graded Lie algebras over  
$^{E}\Omega^{*} \left ({\cal T} \right )$.
The grading is the sum of the
exterior degree and the degree of an $E$-polyvector.
 Fiberwise Schouten bracket also endows   
 $^{E}\Omega^{*} \left ({\cal T}(M) \right )$ with  structure
 of module over the graded Lie algebra 
$^{E}\Omega^{*} \left ( {\cal T} \right )$.
These structures will be respectively denoted by $[,]_{S}$ and ${\displaystyle
  \cdot_{S}}$.
By fiberwise exterior product on $T^{*}_{poly}(\R^{d}_{formal})$,  
$^{E}\Omega^{*} \left ( {\cal T} \right )$ also carries  
a structure of  sheaf of graded commutative algebras
and  $^{E}\Omega^{*} \left ( {\cal T}(M) \right )$ becomes 
a module over the sheaf of graded commutative algebras  
$^{E}\Omega^{*} \left ( {\cal T} \right )$. These structures will both
be denoted by a $\wedge$. 
Thus $^{E}\Omega^{*} \left ( {\cal T}(M) \right )$ is a module
over the  sheaf of
Gerstenhaber algebras $^{E}\Omega \left ( {\cal T} \right )$.

${\bullet}$
Using fiberwise Gerstenhaber bracket, we see that 
$^{E}\Omega^{*} \left ({\cal D} \right )$ is a sheaf of differential
graded Lie algebras and  $^{E}\Omega ({\cal D}(M))$ is a 
module over the sheaf of DGLAs $^{E}\Omega ({\cal D})$. 
These two structures will be denoted $[,]_{G}$ and $\cdot_{G}$.
The grading is the sum of
the exterior degree and the degree of the $E$-polydifferential operator. 
Cuproduct in the space $D^{*}_{poly}(\R ^{d}_{formal})$ endows 
$^{E}\Omega \left ( {\cal D} \right )$ with the structure of 
 a sheaf of DGAAs and 
$^{E} \Omega \left ( {\cal D} (M)\right )$ with the structure
of a module over the sheaf of DGAAs  
$^{E}\Omega \left ( {\cal D} \right )$.\\

$^{E}\Omega ({\cal W})$, $^{E}\Omega ({\cal T}(M))$ and 
$^{E}\Omega ({\cal D}(M))$ are equipped with a decreasing filtration
given  by the order of the
monomials in the fiber coordinates $y^{i}$. \\

In the sequel, we will denote by $\xi^{i}$ the variable $y^{i}$
considered as an element of $\bigwedge ^{1}(E^{*})$. 
Introduce the 2-nilpotent derivation 
$\delta :{}^{E}\Omega ^{*} \left ( {\cal W} \right ) \to 
{}^{E}\Omega ^{*+1}\left ( {\cal W} \right )$ of the sheaf of super
algebras $^{E}\Omega ^{*}\left ( {\cal W} \right )$ defined by 
$\delta = \xi^{i}{\displaystyle \frac{\partial}{\partial y^{i}}}$. 
Using $\cdot _{S}$ and $\cdot_{G}$, $\delta$ extends to a 2-nilpotent
differential of ${\cal T}(M)$ and 
${\cal D} (M)$. 

\begin{proposition}
Let ${\cal B}$ be any of the sheaves  ${\cal W}$, 
${\cal T}(M)$ or ${\cal D}(M)$.
$$H^{\geq 1} \left ( ^{E}\Omega ({\cal B}), \delta \right )=0.$$
Furthermore, we have the following isomorphisms of sheaves of graded
${\cal O}_{X}$-modules~:
$$\begin{array}{l}
H^{0}\left ( ^{E}\Omega ({\cal W}), \delta \right )={\cal O}_{X}\\
H^{0}\left ( ^{E}\Omega ({\cal T}(M)), \delta  \right )=
{}^{E}T_{poly}(M)\\
H^{0}\left ( ^{E}\Omega ({\cal D}^{*}(M)),\delta \right )=
\otimes^{*+1}S(E)
{\displaystyle \mathop \otimes_{{\cal O}_{X}}}M.\\
\end{array}$$
\end{proposition}
 This proposition is  known for ${\cal W}$ and $M={\cal O}_{X}$. 
It is due to Dolgushev 
(~\cite{D1}) for $E=TX$ and to Calaque (~\cite{C1}) for any Lie
 algebroid. Our proof is totally analogous to that of Dolgushev.\\

{\it Proof of the proposition :}

Let us consider the operator $\kappa : {}^{E}\Omega ^{*} ({\cal B}) \to 
{}^{E}\Omega ^{*-1}({\cal B})$ defined by
$$\begin{array}{l}
\forall \sigma \in \Omega^{>0}\left ({\cal T}(M)\right ),\;\;
\kappa (\sigma)= y^{m}{\displaystyle \frac{\partial}{\partial \xi^{m}}}
\int_{0}^{1}\sigma (x, ty, t\xi )
{\displaystyle \frac {dt}{t}}\\
\kappa_{\mid {\cal T}(M)}=0.
\end{array}$$
It satisfies the relation 
$$\delta \kappa + \kappa \delta + {\cal H} =id$$
where 
$$\forall u \in {}^{E}\Omega ^{*}({\cal B}), \;\;\;
{\cal H}(u)=u_{\mid y^{i}= \xi^{i}=0}.$$
The proposition follows.$\Box$\\

{\bf Remark :}

We will keep using the operator $\kappa$ in our proofs. 
Note that $\kappa$ has the two following properties :
\begin{itemize}
\item $\kappa^{2}=0$.
\item $\kappa$ increases the filtration in the variables $y^{i}$'s by one.
\end{itemize}

Let $\nabla$ be a torsion free connection on $E$. 
Let $(e_{1}, \dots, e_{n})$ be a local basis of $E$. Denote by  
$\Gamma_{i,j}^{k}$ the Christoffel symbol of $\nabla$ with respect to
this basis.
As it is explained in previous works
(~\cite{D1}, ~\cite{C1}, ~\cite{D2}, ~\cite{CDH})
such a connection allows to define a
connection  on ${\cal W}$ (still
denoted $\nabla$) as follows :  
$$\nabla = {}^{E}d + \Gamma \cdot \;\;\;
{\rm with}\;\;\;\Gamma = -\xi^{i} \Gamma_{i,j}^{k}y^{j}
{\displaystyle \frac{\partial}{\partial y^{k}}}.$$
It also allows to define a connection on 
${\cal T}(M)$ 
and ${\cal D}(M)$ given by 
$$\nabla _{M}={}^{E}d_{M} + \Gamma \cdot .$$
For example, if 
$\sigma = {\displaystyle \mathop \sum_{l=0}^{\infty}
m_{i_{1},\dots, i_{l}}^{j_{0},\dots,j_{k}}}
y^{i_{1}} \dots y^{i_{l}}
{\displaystyle \frac{\partial}{\partial y^{j_{0}}}}
\wedge \dots \wedge  
{\displaystyle \frac{\partial}{\partial y^{j_{k}}}}$ is a local
section of ${\cal T}(M)$, one has 
$$\begin{array}{rcl}
\nabla_{M}(\sigma )
& = & 
{\displaystyle \mathop {\sum_{l=0}^{\infty}}
{}^{E}d_{M}(m_{i_{1},\dots, i_{l}}^{j_{0},\dots,j_{k}}})
y^{i_{1}} \dots y^{i_{l}}
{\displaystyle \frac{\partial}{\partial y^{j_{0}}}}
\wedge \dots \wedge  
{\displaystyle \frac{\partial}{\partial y^{j_{k}}}}\\
&+& {\displaystyle \mathop \sum_{l=0}^{\infty}
m_{i_{1},\dots, i_{l}}^{j_{0},\dots,j_{k}}}
\Gamma \cdot_{S} \left ( y^{i_{1}} \dots y^{i_{l}}
{\displaystyle \frac{\partial}{\partial y^{j_{0}}}}
\wedge \dots \wedge  
{\displaystyle \frac{\partial}{\partial y^{j_{k}}}}\right ).
\end{array}$$

Since $\nabla$ is torsion free, one has 
$\nabla_{M} \delta + \delta \nabla_{M} =0.$ 
The curvature tensor allows to define the following element of 
$^{E}\Omega^{2}({\cal T}^{0})$ 
$$R=-{\displaystyle \frac {1}{2}}{\xi}^{i}\xi^{j}
\left ( R_{ij} \right )_{k}^{l}(x)y^{k}
{\displaystyle \frac {\partial}{\partial y^{l}}}.$$
A computation shows $\nabla_{M}^{2}=R\cdot :  
{}^{E}\Omega^{*}({\cal B})\to {}^{E}\Omega ^{*+2}({\cal B})$.   
\begin{theorem}
Let ${\cal B}$ be any of the sheaves 
${\cal T}( M)$ and 
 ${\cal D}(M)$. There exists a
 section 
$$A= {\displaystyle \sum_{s=2}^{\infty}} \xi^{k}A^{j}_{k, i_{1}, \dots,
i_{s}}(x)y^{i_{1}}\dots y^{i_{s}}
{\displaystyle \frac {\partial}{\partial y^{j}}}$$
of the sheaf $^{E}\Omega ^{1}({\cal T}^{0})$ such that the
operator  
$D_{M}: ^{E}\Omega ^{*}({\cal B}) \to ^{E}\Omega ^{*+1}({\cal B})$ 
$$D_{M}=\nabla_{M}-\delta + A \cdot $$
is 2-nilpotent and is compatible with the
DG-algebraic structures on $^{E}\Omega ^{*}({\cal B})$.
\end{theorem}
The theorem was proved for ${\cal B}={\cal W}, {\cal T}$ and 
${\cal D}$ in ~\cite{D1} for $E=TX$ 
and in ~\cite{C1} for any algebroid. Our proof is
  inspired by that of ~\cite{D1} (see also ~\cite{C1}).\\

{\it Proof of the theorem}

A computation shows that $D_{M}$ is two-nilpotent if  the
following condition holds :
\begin{equation}
R + \nabla A - \delta A + 
{\displaystyle \frac {1}{2}}[A,A]_{S}=0.
\end{equation}

The following equation 
\begin{equation}
A=\kappa R+\kappa \left ( \nabla (A) +
{\displaystyle \frac {1}{2}}[A,A]_{S}\right )
\end{equation}
has a unique solution (computed by induction on the order  
in the fiber coordinates $y^{i}$'s).
It is shown in ~\cite{D1} 
that the solution of the equation (2) satisfies (1). 
We won't reproduce the proof here.

If $\alpha$ is in $^{E}\Omega ({\cal T})$ and $\mu$ is in 
$^{E}\Omega ({\cal T}(M))$, we have the relations 
$$\begin{array}{l}
D(\alpha \wedge \mu)= D(\alpha)\wedge \mu + 
(-1)^{\mid \alpha \mid +1} \alpha \wedge D(\mu) \\
D(\alpha \cdot_{S}\mu)= D(\alpha)\cdot_{S} \mu + (-1)^{\mid \alpha \mid}
\alpha \cdot_{S} D(\mu)
\end{array}$$
where $\mid \alpha \mid $ denotes the degree of $\alpha$ in the graded
Lie algebra $^{E}\Omega ({\cal T})$. Similarly, 
if $\alpha$ is in $^{E}\Omega ({\cal D})$ and $\mu$ is in 
$^{E}\Omega ({\cal D}(M))$, we have the relations 
$$\begin{array}{l}
D(\alpha \sqcup \mu)= D(\alpha)\sqcup \mu + 
(-1)^{\mid \alpha \mid +1} \alpha \sqcup D(\mu) \\
D(\alpha \cdot_{G}\mu)= D(\alpha)\cdot_{G} \mu + (-1)^{\mid \alpha \mid}
\alpha \cdot _{G} D(\mu)
\end{array}$$
where $\mid \alpha \mid $ denotes the degree of $\alpha$ in the graded
Lie algebra $^{E}\Omega ({\cal D})$.$\Box$\\

One can compute the cohomology of the Fedosov differential $D$.

\begin{theorem}
Let ${\cal B}$ be any of the sheaves 
$^{E}\Omega \left ( {\cal W} \right )$, 
$^{E}\Omega \left ( {\cal T}(M)\right )$ or 
$^{E}\Omega \left ( {\cal D}(M)\right )$,
$$H^{\geq 1}\left ( {\cal B}, D \right )=0.$$
Furthermore, we have the following isomorphisms of sheaves of graded
commutative algebras 
$$\begin{array}{l}
H^{0} \left ( ^{E}\Omega ({\cal W}), D\right )\simeq{\cal O}_{X}\\
H^{0}\left ( ^{E}\Omega ({\cal T}), D\right )\simeq \bigwedge^{*+1}E\\
\end{array}$$
and the following isomorphism of sheaves of DGAAs (over $\R$)
$$H^{0}\left ( ^{E}\Omega ({\cal D}), D\right )\simeq 
\otimes^{*+1}S(E).$$
Using the identification above, 
$H^{0}\left ( ^{E}\Omega ({\cal T}(M)), D\right )$ and 
$\bigwedge^{*+1}E {\displaystyle \mathop \otimes_{{\cal O}_{X}}} M$
are isomorphic as 
$H^{0}\left ( ^{E}\Omega ({\cal T}), D\right )\simeq
\bigwedge^{*+1}E$-modules. Furthermore 
$H^{0}\left ( ^{E}\Omega ({\cal D}(M)), D\right )$
and $\otimes^{*+1}S(E)
{\displaystyle \mathop \otimes_{{\cal O}_{X}}}M$ are isomorphic
as $H^{0}\left ( ^{E}\Omega ({\cal D}), D\right )\simeq 
\otimes^{*+1}S(E)$-modules.
\end{theorem}
This theorem is already known for $M={\cal O}_{X}$: see ~\cite{D1} for
the case where $E=TX$  and ~\cite{C1}, ~\cite{C2} 
for any Lie algebroid. The proof
of the theorem is very similar to the proof in the case where 
$ M={\cal O}_{X}$. That is why, we give only a sketch of it 
 and refer to ~\cite{CDH} and ~\cite{C2} for details.\\

{\it Proof of the theorem :}

The first assertion of the theorem follows from a spectral sequence
argument using the filtration  
on ${\cal B}$ given by the order on the $y^{i}$'s 
(see ~\cite{CDH}, theorem 2.4 for details).

Let $u \in  {\cal B}  \bigcap Ker \delta$. 
One can show (solving the equation by induction on the order in the
fiber coordinates $y^{i}$'s) 
that there exists 
a unique $\lambda (u) \in {\cal B} \bigcap Ker D$ such that 
$$\lambda (u)=u + \kappa \left ( \nabla \lambda (u) + 
A \cdot \lambda (u)\right ).$$
 Thus, we have
defined a map 
$\lambda : Ker \delta \bigcap {\cal B} \to Ker D \bigcap {\cal
  B}$. One can show that $\lambda$ is bijective and that 
$\lambda ^{-1}={\cal H}$. The following relations (easy to establish)
allows to finish the proof of the theorem :\\

${\bullet}$ If $\alpha , \beta \in {}^{E}\Omega ({\cal W})$, then 
${\cal H}(\alpha \beta)={\cal H}(\alpha){\cal H}(\beta)$\\

$\bullet$ If $\alpha \in {}^{E}\Omega ({\cal T})$ and 
$\mu \in {}^{E}\Omega ({\cal T}(M))$, then 
${\cal H} (\alpha \wedge \mu) ={\cal H}(\alpha) \wedge {\cal H}(\mu)$.\\

$\bullet$ If $\alpha \in {}^{E}\Omega ({\cal D})$ and 
$\mu \in {}^{E}\Omega ({\cal D}(M))$, then 
${\cal H} (\alpha \sqcup \mu) ={\cal H}(\alpha) \sqcup {\cal H}(\mu).$
$\Box$\\

As $D$ is compatible with the action 
$\cdot_{S}$ of $^{E}\Omega^{*}\left ({\cal T} \right )$ over 
${}^{E}\Omega ^{*} \left ({\cal T}(M) \right )$ and hence with
the Schouten bracket on 
${}^{E}\Omega^{*}\left ({\cal T} \right )$, 
$H^{*}\left ( ^{E}\Omega ({\cal T}), D \right )$ is a
graded Lie algebra and 
$H^{*}\left ( ^{E}\Omega ({\cal T}(M)), D \right )$
is a module over the graded Lie algebra 
$H^{*}\left ( ^{E}\Omega^{*} ({\cal T}), D \right )$. 
So, it is  natural to wounder whether the isomorphisms of the previous
proposition respect this structure.   
\begin{proposition}
The map ${\cal H} : {\cal T}^{*} \bigcap Ker D \to 
{\cal T}^{*}\bigcap Ker \delta \simeq {}^{E}T^{*}_{poly}$ is an
isomorphism  of graded Lie algebras.

The map ${\cal H} : {\cal T}^{*}(M) \bigcap Ker D \to 
 {}^{E}T^{*}_{poly}(M)$ is an
isomorphism of modules over the graded Lie algebras 
${\cal T}^{*} \bigcap Ker D \simeq {}^{E}T^{*}_{poly}$.
\end{proposition}

{\it Proof of the proposition :}

The first assertion of the proposition is proved in 
~\cite{C1}, ~\cite{C2}. 
 Let us now prove the second assertion. 
Denote by $\pi$ the map  from $D(E)$ to $End(M)$ defined by the action of
$D(E)$ on $M$. 

Let $m$ be an element of $M$ and let 
$u={\displaystyle \mathop \sum_{i=1}^{d}u_{i}(x)e_{i}} \in 
{}^{E}T^{0}_{poly}$. 
Using the definition of 
$\lambda$, one finds easily :
$$\begin{array}{l}
\lambda (m) = m +
{\displaystyle \sum_{i=1}^{d}}
y^{i}\pi (e_{i})\cdot m \;\;{\rm mod}\;\; \mid y \mid  \\
\lambda (u)= {\displaystyle \mathop \sum_{i=1}^{d}}u_{i}
{\displaystyle \frac {\partial}{\partial y^{i}}} \;\;{\rm mod}\;\; 
\mid y \mid \\
\end{array}$$
Hence $$\lambda (u) \cdot \lambda (m) = 
{\displaystyle \mathop \sum_{i=1}^{d}}u_{i}\pi(e_{i})\cdot m 
\;\;{\rm mod}\;\; \mid y \mid .$$
and 
$${\cal H}(\lambda (u) \cdot \lambda (m))= u \cdot m =
{\cal H}(\lambda (u)) \cdot {\cal H}(\lambda (m)).$$
The end of the proof follows from the definition of the action of 
${}^{E}T_{poly}$ on ${}^{E}T_{poly}(M)$ and the previous
theorem. $\Box$\\

{\bf The morphism $\mu '_{M}$}

Let us first recall the construction of $\mu '$ (~\cite{CDH}). 
${\cal T}^{0}$ is the sheaf of Lie algebras over the sheaf of algebras
${\cal T}^{-1}=\widehat{S}(E^{*})$ and we have 
${\cal D}^{0}= D({\cal T}^{0})$. 
 The morphism of Lie algebras 
$\lambda = {\cal H}^{-1} : E \to {\cal T}^{0} \bigcap Ker D $ induces
a morphisms of sheaves of algebras $\mu : D(E) \to {\cal D}^{0}$ that takes
values in $Ker D \bigcap {\cal D}^{0}$. 
We will denote by  $\mu '$ the only morphism of sheaves of DGAAs 
from  ${}^{E}D_{poly}^{*}$ to ${\cal D}^{*} $ defined by 
$$\mu ' _{\mid {}^{E}D^{0}_{poly}}= \mu ,\;\;\;\;
\mu ' _{\mid {\cal O}_{X}}=\lambda .$$
Let $\mu '_{M} : {}^{E}D_{poly}^{*}(M) \to 
{\cal D}^{*}(M)$  
the morphism defined by : $ \forall P_{0}, \dots, P_{n} \in D(E), \;\;
\forall m \in M$
$$\begin{array}{l}
{\mu '_{M}}(m)=\lambda (m) \\
\mu '_{M}(P_{0} \otimes \dots \otimes P_{n} \otimes m) =
\mu (P_{0}) \otimes \dots \otimes \mu(P_{n}) \otimes \lambda (m)\\
\end{array}$$
Note that $\mu ' =\mu '_{{\cal O}_{X}}$.
\begin{proposition}
a) $\mu$ is an isomorphism of sheaves of algebras from $D(E)$  to 
${\cal D}^{0} \bigcap Ker D$. It is also a morphism of sheaves of
bialgebroids. 

b) $\mu '$ is an isomorphism of sheaves of DGLAs from ${}^{E}D_{poly}^{*}$ to 
${\cal D}^{*} \bigcap Ker D$. It is also an isomorphism of sheaves of
DGAAs.

c) $\mu '_{ M} : 
 {}^{E}D_{poly}^{*}\left (M \right ) \to 
{\cal D}^{*}\left ( M \right ) \bigcap Ker D$ is an isomorphism
 of modules over the sheaf of DGLAs  
${}^{E}D_{poly}^{*} \simeq {\cal D}^{*}\bigcap Ker D$. 
It is also an isomorphism of modules over the sheaf of 
DGAAs ${}^{E}D_{poly}^{*} \simeq {\cal D}^{*}\bigcap Ker D$.
\end{proposition}

{\it Proof of the proposition :}

a) and b) are shown in ~\cite{CDH}.
The proof of c) is analogous. 
Using the definition of $\mu $ and $\mu'_{M}$, one can show
easily the following :
$$\forall P \in D(E),\; \forall m \in M ,\;\; 
\mu '_{M}(P \cdot m)= \mu (P)\cdot \mu'_{M}(m).$$
As moreover $\mu $ is an isomorphism of bialgebroids (~\cite{CDH}), 
$\mu '_{M}$ is a morphism of modules over the sheaf of DGLAs  
${}^{E}D_{poly}^{*} \simeq {\cal D}^{*}\bigcap Ker D$.
$\mu '_{M}$ is clearly a morphism
of modules over the sheaf of DGAAs  
${}^{E}D_{poly}^{*} \simeq {\cal D}^{*}\bigcap Ker D$.
 The fact that  $\mu ' _{M}$ an isomorphism
follows from a) and theorem 4.1.3. $\Box$

\subsection{Kontsevitch's result}

Recall that $\R ^{d}_{formal}$ is the formal completion of 
$\R ^{d}$ at the
origin. The ring of functions on $\R^{d}_{formal}$ is 
$\R [[y^{1}, \dots, y^{d}]]$ and the Lie-Rinehart algebra of vector
fields is $Der \left ( \R [[y^{1}, \dots, y^{d}]] \right )$. Denote by 
$T_{poly}^{*}(\R ^{d}_{formal})$ and 
$D_{poly}^{*}(\R^{d}_{formal})$ the 
DGLAs of polyvector fields and polydifferential
operators on  $\R ^{d}_{formal}$ respectively. 

\begin{theorem}\label{Kontsevitch}
There exists a quasi-isomorphism of $L_{\infty}$ algebras $U$ from 
$T^{*}_{poly}(\R ^{d}_{formal})$ to 
$D^{*}_{poly}(\R^{d}_{formal})$ 
such that 

(1) The first structure map $U^{[1]}$ is the quasi-isomorphism
    $U_{HKR}$.

(2) $U$is $GL_{d}(\R)$-equivariant.

(3) If $n>1$ then, for any vector fields 
$v_{1}, \dots, v_{n} \in T^{0}_{poly}(\R ^{d}_{formal})$
$$U^{[n]}(v_{1}, \dots, v_{n})=0$$

(4) If $n>1$ then  for any vector field $v$ linear in the coordinates $y^{i}$ 
and polyvector fields 
$\chi_{2}, \dots, \chi_{n} \in T^{*}(\R^{d}_{formal})$
$$U^{[n]}(v, \chi_{2}, \dots, \chi_{n})=0.$$
\end{theorem}
Moreover, Kontsevitch gives an explicit expression for $U^{[n]}$ 
(~\cite{Ko}, see also ~\cite{AMM} or ~\cite{BCKT} 
for a detailed exposition) 
which involves admissible graphs.

\begin{definition}
Let $n$ and $m$ be two integers. 
An admissible graph $\Gamma$ of type $(n,m)$ is a labeled oriented 
graph satisfying the following properties. Let $V_{\Gamma}$ be the
set of vertices of $\Gamma$ and $E_{\Gamma}$ be the set of edges of 
$\Gamma$.

1)  $V_{\Gamma} = \{1, \dots,n \} \sqcup 
\{\bar{1}, \dots , \bar{m} \}$. Elements of  $\{1, \dots, n\}$ are
called first type vertices and element of 
$\{\bar{1}, \dots, \bar{m} \}$ second type  vertices.    

2) Every edge of $\Gamma$ starts from a first type vertex.

3) There is no loop.

4) Two edges can't have the same source and the same target.
\end{definition} 

We will write $G_{n,m}$ for the set of admissible graphs with $n$
first type vertices and $m$ second type vertices. Let $\Gamma$ be an
element of $G_{n,m}$. We will denote by $E_{\Gamma}$ the set of its
edges. If $\gamma$ is in $E_{\Gamma}$, then $s(\gamma)$ will be its
source and $t(\gamma)$ its target. 
Let us introduce the following notation : if $k$ is a vertex of
first type 
$$\begin{array}{l}
(k,*)=\{\gamma  \in E_{\Gamma} \mid s(\gamma)=k \}=
\{e_{k}^{1},\dots, e_{k}^{s_{k}} \}  \\
\end{array}$$
Similarly, the subset $(*,k)$ of $E_{\Gamma}$ is defined for any
vertex of $\Gamma.$ 

Let $\alpha_{1}, \dots, \alpha_{n}$ be $n$ polyvector fields such that for
any $j \in [1,n]$, $\alpha_{j}$ is a $s_{j}$ polyvector fields.
We will associate to such $\alpha_{1}, \dots, \alpha_{n}$ an 
$m$ polydifferential operator $B_{\Gamma}(\alpha_{1}, \dots, \alpha_{n})$.
Write 
$$\alpha_{j}=
{\displaystyle \mathop \sum_{i_{1},\dots,i_{s_{j}}}}
\alpha^{i_{1}, \dots, i_{s_{j}}}\partial_{i_{1}} \wedge \dots
\wedge \partial_{i_{s_{j}}} \; 
{\rm with} \; 
\partial_{k}=
{\displaystyle \frac{\partial}{\partial y^{k}}}.$$
If $I :E_{\Gamma}\to \{1, \dots, d\}$ is a map from $E_{\Gamma}$ to
$\{1. \dots, d\}$, we set 
$$\begin{array}{l}
D_{I(x)}={\displaystyle \prod_{e \in (*,x)}}\partial_{I(e)}\\
\alpha_{k}^{I}=\alpha_{k}^{I(e_{k}^{1}), \dots, I(e_{k}^{s_{k}})}.
\end{array}$$
$B_{\Gamma}(\alpha_{1} \otimes  \dots \otimes \alpha_{n})$ 
is the $m$-differential operator
defined\  by : for any functions $f_{1}, \dots, f_{m}$, 
$$B_{\Gamma}
\left (\alpha_{1} \otimes  \dots \otimes \alpha_{n}  \right )
(f_{1}, \dots, f_{m})=
{\displaystyle \sum_{I:E_{\Gamma}\to \{1, \dots, d\}}}
{\displaystyle \prod_{k=1}^{k=n}}D_{I(k)}\alpha_{k}^{I}
{\displaystyle \prod_{l=1}^{l=m}}D_{I(\bar{l})}f_{l}$$ 

If $\alpha_{1}, \dots, \alpha_{n}$ are any graded elements of
$T_{poly}$, one has 
$$U^{[n]}(\alpha_{1}, \dots, \alpha_{n})=
{\displaystyle \mathop \sum_{\Gamma \in G_{n,m}}}
W_{\Gamma}B_{\Gamma}(\alpha_{1}\otimes \dots \otimes \alpha_{n})$$
where the sum is taken over the graph $\Gamma$ in $G_{n,m}$ such that 
$B_{\Gamma}(\alpha_{1}\otimes \dots \otimes \alpha_{n})$ is defined
and the relation $m-2+2n={\displaystyle \sum_{i=1}^{n}s_{k}}$ is
satisfied. 
The coefficient $W_{\Gamma}$ can be different from zero
only if $\mid E_{\Gamma} \mid =2n + m -2.$ Let us now describe
it.

Let ${\cal H}$ be the Poincar{\'e} half plane :
$${\cal H}=\{z \in \C \mid Im(z)>0\}.$$ 
Introduce 
$$Conf_{n,m}=\{(p_{1}, \dots, p_{n}, q_{\bar{1}}, \dots, q_{\bar{m}}) \in {\cal
  H}^{n}\times \R^{m} \mid 
 p_{i} \neq p_{j}, \;  q_{\bar{i}}\neq q_{\bar{j}}\}.$$ 
The group 
$G=\{z \mapsto az+b \mid (a,b) \in \R^{+*} \times \R\}$
acts freely on $Conf_{n,m}$. The quotient $C_{n,m}=Conf_{n,m}/G$ is a
manifold of dimension $2n+m-2$. As $Conf_{n,m}$ is naturally oriented
and the action of $G$ preserves this orientation, $C_{n,m}$ inherits a
natural orientation. 
$C_{n,m}$ has several connected components, we will use one of them
$C_{n,m}^{+}$ defined by 
$$C_{n,m}^{+}=\{(p_{1}, \dots, p_{n},q_{\bar{1}}, \dots,q_{\bar{m}}) \mid 
q_{\bar{1}} <\dots <q_{\bar{m}} \}.$$
If $i\in \{1, \dots,n \}$ and 
$j \in \{1,\dots,n\} \sqcup \{\bar{1},\dots, \bar{m}\} $ 
(with $i \neq j$),
one defines a function 
$$\begin{array}{rcl}
\theta_{i,j} : C_{n,m} & \to & \R/2\pi \Z\\
              (z_{k})_{k\in[1,n] \sqcup [1,\bar{m}]} & \mapsto & 
{\displaystyle \frac{1}{2 \pi}}Arg
{\displaystyle \frac{z_{j}-z_{i}}{z_{j}- \overline{z_{i}}}}
\end{array}.$$
Let   $\Gamma$  be an element of $G_{n,m}$. We order $E_{\Gamma}$ with
the lexicographic order and define  the closed form 
$$ \omega_{\Gamma}={\displaystyle \mathop \wedge_{\gamma \in E_\Gamma}}
d\theta_{s(\gamma), t(\gamma)}.$$ One then put 
$$W_{\Gamma}=\int_{C_{n,m}^{+}}\omega_{\Gamma}.$$
This integral is absolutely convergent as the integrand  extends to a
differential form on a 
compactification of $C_{n,m}^{+}$, $\overline{C_{n,m}^{+}}$, which is
a manifold with corners of dimension $2n +m-2$
 (~\cite{Ko}, see also ~\cite{AMM} and ~\cite{BCKT}). 

\begin{lemma}~\label{ouvert de coordonnees}
Let $n$ be a non zero integer. 
For any polyvectorfields $\gamma_{1}, \dots, \gamma_{n}$, one has 
$$U^{[n+1]}({\displaystyle \frac{\partial}{\partial y^{i}}}, 
\gamma_{1}, \dots, \gamma_{n})=0.$$
\end{lemma}

{\it Proof of the lemma :}

We will prove that for any $\Gamma$ in $G_{n+1,m}$ having a
contribution in $U^{[n+1]}$, one has $W_{\Gamma}=0$. For such a
$\Gamma$, there is no edge going to the vertex $1$ and there is
exactly one edge starting from the vertex $1$ and going to a vertex
$i_{0}$ which might be of first or of second type. We will denote by
$\Gamma '$  the element of $G_{n,m}$ obtained from $\Gamma$ by removing 
the vertex $1$ and the edge going from $1$ to $i_{0}$.\\ 

{\it First case : $i_{0}$ is of first type }

Using the action of $G$, we put $p_{i_{0}}$ in $i$. If $j$ is in 
$[1,n+1]-\{i_{0}\}$, we will write  
$z_{j}=a_{j} +ib_{j}$ for the affix of $p_{j}$ and if $k$ is in
$[1,m]$, we will write $t_{k}$ for the coordinate of $q_{k}$. 
One has $$\omega_{\Gamma}=
{\displaystyle \frac{1}{2\pi}}
d Arg \left ({\displaystyle \frac{i-z_{1}}{i-\overline{z_{1}}}} \right
) \wedge \omega_{\Gamma '}$$
and $\omega_{\Gamma '}$ is a differential form of degree $2(n+1)+m-3$
 in the $2(n-1)+m$ variables $da_{2}, db_{2}, \dots,  
\widehat{da_{i_{0}}}, 
\widehat{db_{i_{0}}}, \dots, da_{n+1}, db_{n+1}, dt_{1}, \dots,
dt_{m}$. Hence $\omega_{\Gamma '}=0$ and $\omega_{\Gamma}=0$. \\

{\it Second case : $i_{0}$ is of second type }

We treat the case where $i_{0}\neq \bar{m}$. The case where $i_{0}
=\bar{m}$ is treated analogously.  
Using the action of $G$, we put $q_{i_{0}}$ in $0$ and $q_{i_{0}+1}$
in $1$. One has 
$$\omega_{\Gamma}=
{\displaystyle \frac{1}{\pi}}
dArg(z_{1})\wedge \omega_{\Gamma '}.$$
$\omega_{\Gamma '}$ is a differential form of degree $2(n+1)+m-3$ in the 
$2n+m-2$ variables $a_{2},b_{2}, \dots, a_{n+1}, b_{n+1}, q_{1},
\dots, \widehat{q_{i_{0}}},\widehat{q_{i_{0}+1}}, \dots, q_{m}$. 
Hence $\omega_{\Gamma '}=0$ and $\omega_{\Gamma}=0$. $\Box$
 
\subsection{Proof of the formality theorem}

The proof will follow ~\cite{C2}. Before starting the proof, 
let's recall the following well known fact of sheaf theory : 
If ${\cal C}^{*}_{1}$ and ${\cal C}_{2}^{*}$ are complexes of c-soft
sheaves and if $\Theta$ is a quasi-isomorphism from ${\cal C}_{1}^{*}$
to ${\cal C}_{2}^{*}$, then $\Gamma \left (\Theta \right )$ is a
quasi-isomorphism from $\Gamma \left ({\cal C}_{1}^{*}\right )$ to 
$\Gamma \left ({\cal C}_{2}^{*}\right )$.\\

We will adopt the following notations :

$\lambda_{T}^{M} : {}^{E}T^{*}_{poly}(M) \to 
{}^{E}\Omega \left ( {\cal T}(M)\right )$ is the inverse of the
map ${\cal H}$.

$\lambda_{D}^{M} : {}^{E}D^{*}_{poly}(M) \to 
{}^{E}\Omega \left ( {\cal D}(M)\right )$ is the map
$\mu'_{M}$. 

We set $\lambda_{D}^{{\cal O}_{X}}=\lambda _{D}$ and 
$\lambda_{T}^{{\cal O}_{X}}=\lambda _{T}$. From Kontsevitch's work
(theorem ~\ref{Kontsevitch}), we
know that there exists a fiberwise quasi-isomorphism of $L_{\infty}$-algebras
${\cal U}$  
from $^{E}\Omega ({\cal T})$ to $^{E}\Omega ({\cal D})$ 
whose Taylor coefficients
will be denoted 
${\cal U}^{[n]}: S^{n}\left (^{E}\Omega ({\cal T})[1] \right ) \to 
^{E}\Omega ({\cal D})$ (first we construct ${\cal U}$ on an open
subset trivializing $E$ and then glue the $L_{\infty}$-morphisms). 
Using the explicit expression  of $U^{[n]}$ 
(~\cite{Ko}, ~\cite{AMM}), 
one sees easily
that ${\cal U}^{[n]}$ still make sense if we replace the last
argument by an element of $^{E}\Omega \left ({\cal T}(M)\right )$. 
Thus we define 
${\cal V}^{[n]} : 
S^{n} \left (^{E}\Omega ({\cal T})[1] \right ) \otimes 
^{E}\Omega \left ({\cal T}(M) \right ) 
\to ^{E}\Omega \left ( {\cal D}(M)\right )$ by 
$$\forall \gamma_{1}, \dots, \gamma_{n} \in ^{E}\Omega ({\cal T})[1], 
\forall \nu \in ^{E}\Omega \left ({\cal T}(M) \right ),\;\;\;
{\cal V}^{[n]}(\gamma_{1}, \dots, \gamma_{n}, \nu)=
{\cal U}^{[n+1]} (\gamma_{1}, \dots, \gamma_{n}, \nu).$$

Thus we get the following diagram 


$$\xymatrix{
\left ( {}^{E}\Omega \left ( {\cal T} \right ), 0, [,]_{S} 
\right ) 
\ar[d]_{\cdot_{S}}^{L_{\infty}-mod}
\ar[r]^{\cal U}
& \left ( {}^{E}\Omega \left (  {\cal D} \right ),\partial, [,]_{G} 
\right ) 
\ar[d]_{\cdot_{G}}^{L_{\infty}-mod}\\
\left ( {}^{E}\Omega \left (  {\cal T}(M) \right ), 0, \cdot_{S} 
\right )  \ar[r]^{\cal V}
&\left ( {}^{E}\Omega \left (  {\cal D}(M)\right ), \partial_{M} , \cdot_{G} 
\right ) 
}$$

Let $V$ be an open subset on which $E_{\mid V}$ is trivial. The
differential $^{E}d$ (respectively $^{E}d_{M}$) is defined on 
${}^{E}\Omega \left (  {\cal T} \right )_{\mid V}$ and  
${}^{E}\Omega \left (  {\cal D} \right )_{\mid V}$ 
(respectively ${}^{E}\Omega \left (  {\cal T}(M) \right )_{\mid V}$
and 
${}^{E}\Omega \left (  {\cal D}(M) \right )_{\mid V}$). 
As the quasi-isomorphisms of the previous diagram are fiberwise, we
can add the differentials $^{E}d$ and $^{E}d_{M}$, in the previous
quasi-isomorphism. We get a morphism of $L_{\infty}$-algebras 
$$\overline{\cal U}: \;\;
\left ( {}^{E}\Omega \left (  {\cal T} \right )_{\mid V}, ^{E}d, [,]_{S} 
\right ) \to
 \left ( {}^{E}\Omega 
\left (  {\cal D} \right )_{\mid V}, ^{E}d + \partial , [,]_{G}  \right )$$  
and a morphism of $L_{\infty}$-modules 
over 
$^{E} \Omega \left ({\cal T}\right )_{\mid V} $ 
$$\overline{\cal V} : \;\; 
\left ( {}^{E}\Omega \left (  {\cal T}(M)\right )_{\mid V}, 
^{E}d_{M}, \cdot_{S}  \right )\to
\left ( {}^{E}\Omega \left ({\cal D}(M)\right ) _{\mid V}, 
^{E}d_{M}+ \partial_{M} , \cdot_{G} \right )$$ 
We endow  ${\cal B}={\cal T}(M)_{\mid V}$ or ${\cal D}(M)_{\mid V}$ 
with the  filtration 
$$F^{p}\left (^{E}\Omega ({\cal B}) \right )=
{\displaystyle \mathop \oplus _{k \geq p}}^{E}\Omega ^{k}({\cal B})$$
A spectral sequence argument shows that 
$\overline{\cal U}$ and $\overline{\cal V}$ are quasi-isomorphisms
(see ~\cite{C2} and ~\cite{CDH} for details).
Thus, we have the following diagram where the horizontal arrows are
quasi-isomorphisms 
 $$\xymatrix{
\left ( {}^{E}\Omega \left (  {\cal T} \right )_{\mid V}, ^{E}d, [,]_{S} 
\right ) 
\ar[d]_{\cdot_{S}}^{L_{\infty}-mod}
\ar[r]^{\overline{\cal U}}
& \left ( {}^{E}\Omega \left (  {\cal D} \right )_{\mid V}, 
^{E}d + \partial , [,]_{G} 
\right ) 
\ar[d]_{\cdot_{G}}^{L_{\infty}-mod}\\
\left ( {}^{E}\Omega \left (  {\cal T}(M)\right )_{\mid V} , 
^{E}d_{M}, \cdot_{S} 
\right )  \ar[r]^{\overline{\cal V}}
&\left ( {}^{E}\Omega \left (  {\cal D}(M)\right )_{\mid V}, 
^{E}d_{M} + \partial_{M} , \cdot_{G} \right ) 
}$$

On  $V$, the Fedosov
differential can be written $D_{M}={}^{E}d_{M}+B$ with 
$$ B= {\displaystyle \mathop \sum_{p=0}^{\infty}}\xi^{i}
B_{i,j_{1}, \dots, j_{p}}(x) 
y^{j_{1}}\dots y^{j_{p}}
{\displaystyle \frac{\partial}{\partial y^{k}}}.$$
We set $D=D_{{\cal O}_{X}}$. The element 
$B$ of  $^{E}\Omega^{1} ({\cal T}^{0})_{\mid V}$ is a Maurer Cartan element
of the (filtered) sheaf of DGLAs 
$\left ( ^{E}\Omega({\cal T})_{\mid V}, ^{E}d, [,]_{S} \right )$. 
This means that 
$\left ( ^{E}\Omega ({\cal T}( M))_{\mid V}, D_{M}, \cdot_{S} \right
)$ is obtained from 
$\left ( ^{E}\Omega ({\cal T}(M))_{\mid V}, ^{E}d_{M}, \cdot_{S} \right
)$ via the twisting procedure by the Maurer Cartan element $B$ 
(~\cite{D2}).
We know that ${\displaystyle \sum_{n \geq 1}\frac{{\cal U}^{[n]}(B^{n})}{n!}}$ 
is a Maurer Cartan section of 
$\left ( ^{E}\Omega ({\cal D})_{\mid V}, {}^{E}d + \partial ,
  \cdot_{G} \right )$. But, due to  property (3) of $U$, 
${\displaystyle \sum_{n \geq 1}\frac{{\cal U}^{[n]}(B^{n})}{n!}}=B$. 
Twisting $\overline{\cal U}$ and $\overline{\cal V}$ by the Maurer Cartan
element $B$ (\cite{D2}), 
we get the following diagram where the horizontal arrows
are quasi-isomorphism

 $$\xymatrix{
\left ( {}^{E}\Omega \left ( {\cal T} \right )_{\mid V}, D, [,]_{S} 
\right ) 
\ar[d]_{\cdot_{S}}^{L_{\infty}-mod}
\ar[r]^{\overline{\cal U}^{B}}
& \left ( {}^{E}\Omega \left ( {\cal D} \right )_{\mid V}, 
D+ \partial , [,]_{G} 
\right ) 
\ar[d]_{\cdot_{G}}^{L_{\infty}-mod}\\
\left ( {}^{E}\Omega \left ({\cal T}(M)\right )_{\mid V}, 
D_{M}, \cdot_{S} 
\right )  \ar[r]^{\overline{\cal V}^{B}}
&\left ( {}^{E}\Omega \left ({\cal D}(M)\right )_{\mid V}, 
D_{M}+\partial_{M} , \cdot_{G} \right ) 
}$$
$\overline{\cal U}^{B}$ and $\overline{\cal V}^{B}$ do not depend on
the choice of the trivialization of $E_{\mid V}$ and hence are a well
defined morphisms of $L_{\infty}$-algebras and $L_{\infty}$-modules 
respectively. Indeed the only term in $B$ that depends on the
coordinates is 
$\Gamma = - \xi^{i}\Gamma_{i,j}^{k}y^{j}
{\displaystyle \frac{\partial}{\partial y^{k}}}$
and it is linear in the fiber coordinates $y^{i}$ so that it does
neither contribute to $\overline{\cal U}^{B}$ nor to 
$\overline{\cal V}^{B}$ thanks to property (4) of $U$ 
(see ~\cite{D1},~\cite{C1},~\cite{D2},~\cite{CDH} for details). 
Hence $\overline{\cal U}^{B}$ and 
$\overline{\cal V}^{B}$ are defined globally and we get the following diagram.
 
 $$\xymatrix{
\left ( {}^{E}\Omega \left ( {\cal T} \right ), D, [,]_{S} 
\right ) 
\ar[d]_{\cdot_{S}}^{L_{\infty}-mod}
\ar[r]^{\overline{\cal U}^{B}}
& \left ( {}^{E}\Omega \left ( {\cal D} \right ), 
D+ \partial , [,]_{G} 
\right ) 
\ar[d]_{\cdot_{G}}^{L_{\infty}-mod}\\
\left ( {}^{E}\Omega \left ({\cal T}(M) \right ), 
D_{M}, \cdot_{S} 
\right )  \ar[r]^{\overline{\cal V}^{B}}
&\left ( {}^{E}\Omega \left ({\cal D}(M)\right ), 
D_{M}+\partial_{M} , \cdot_{G} \right ) 
}$$

The following lemma shows that the map
$\lambda_{D}^{M}(X)$ (and hence $\lambda_{D}(X)$) is a quasi-isomorphism
 from 
$\left [ \Gamma \left (^{E}D_{poly}(M)\right ), \partial_{M} \right ]$ 
to 
$\left  [\Gamma \left (^{E}\Omega\left ({\cal D}(M) \right )\right ), 
D_{M}+ \partial_{M} \right ]$.

\begin{lemma}
 The natural inclusion $\iota :$
$\left [ \Gamma \left ({\cal D}^{*}(M) \bigcap Ker D_{M} \right ) 
, \partial_{M} \right ]
\hookrightarrow \linebreak  
\left [ \Gamma \left (\Omega ^{*} ({\cal D}(M))\right ) , 
D_{M} + \partial_{M} \right ]$
is  a quasi-isomorphism. 

\end{lemma}

{\it Proof of the lemma :}

Consider a decomposition of $Ker (D_{M} + \partial_{M})$ of the form 
$$Y \oplus  Im (D_{M}+\partial_{M}) = Ker (D_{M} + \partial_{M}).$$
One may construct a map 
$V : Ker (D_{M}+ \partial_{M} ) \to 
\Gamma \left ( \Omega \left ({\cal D}(M) \right )\right )$ such that 

i) for any $x$ in $Ker (D_{M}+\partial_{M} )$, 
$x - (D_{M} + \partial_{M} )(V(x)) \in 
\Gamma \left ({\cal D}(M)\bigcap Ker D_{M} \right )$

ii) If $x \in Im (D_{M} + \partial _{M} )$, $V(x)$ is a preimage of $x$ by 
$D_{M}+ \partial_{M}$. 

It is enough to construct $V(x)$ for $x$ in $Y$. Write 
$x=x_{r}+ \dots + x_{0}$ with 
$x_{i} \in \Gamma \left (\Omega^{i} \left ( {\cal D}(M) \right )\right )$. 
The equality $(D_{M}+ \partial_{M} )(x)=0$ implies $D_{M}(x_{r})=0$
(because $\partial_{M}$ preserves the exterior degree). 
Then using
the exactness of $D_{M}$, we construct a map 
$V_{r} : Y \to 
\Gamma \left (\Omega^{\leq r-1}\left ({\cal D}(M) \right )\right )$ 
such that for any $x$ in $Y$, 
$x - (D_{M} + \partial_{M})V_{r}(x)$ has maximal exterior degree inferior
or equal to $r-1$. Going on like this, we construct $V$. 

We may now exhibit an inverse to $H^{i}(\iota)$. With obvious
notations, we have  
$$H^{i}(\iota)^{-1} \left ( [ \mu] \right )= [ \mu - (D_{M} + \partial_{M}
)V(\mu) ].$$

This finishes the proof of the lemma.$\Box$ \\

As $\lambda_{D}^{M}(X)$ is a quasi-isomorphism of $L_{\infty}$-modules
over $\Gamma \left (^{E}D^{*}_{poly}\right )$, 
 there \linebreak
exists a quasi-isomorphism of $L_{\infty}$-modules
over $\Gamma \left (^{E}D^{*}_{poly}\right )$,\linebreak  
$\alpha_{D}^{M} : 
\left [ \Gamma \left ({}^{E}\Omega \left ( {\cal D}\left ( M \right )\right
  )\right ), D_{M}+\partial_{M}\right ] 
\to 
\left [ \Gamma \left (^{E}D_{poly}^{*}(M)\right ), \partial_{M}\right ]$ 
such that \linebreak 
$H^{i}\left ( \alpha ^{M[1]}_{D} \right )= 
H^{i} \left ( \lambda_{D}^{M}  \right )^{-1}$ 
(see ~\cite{AMM} for the case of
$L_{\infty}$ algebras). 
The morphism 
${\cal V}_{M}=\alpha _{D}^{M} \circ \overline{\cal V}^{B}(X) \circ 
\lambda_{T}^{M}(X)$ is a quasi-isomorphism of $L_{\infty}$-modules
 over $\Gamma \left (^{E}T^{*}_{poly} \right )$ 
from $\Gamma \left ({}^{E}T^{*}_{poly}(M)\right )$ 
to $\Gamma \left ( {}^{E}D^{*}_{poly}(M)\right )$. 
One checks easily that ${\cal V}_{M}^{[0]}$ induces $U^{M}_{HKR}$ in
cohomology.

Inverting  $\lambda_{D}$ into a quasi-isomorphism of $L_{\infty}$
algebras provides Calaque 's quasi-isomorphism of $L_{\infty}$
algebras ${\Upsilon}$ from 
$\Gamma \left (^{E}T^{*}_{poly}\right )$ to 
$\Gamma \left (^{E}D^{*}_{poly}\right )$ (~\cite{C2}). This finishes the
proof of the theorem \ref{formality}. $\Box$

\subsection{Local expression of ${\cal V}_{M}$ in the case of the
  tangent bundle of $\R^{d}$.}

In this section, assume that $X=\R^{d}$ and $E=T\R^{d}$. 
We choose the connection whose Christoffel symbols are
$0$. Thus, we have 
$$\nabla (f {\displaystyle \frac{\partial}{\partial x^{i}}})=
df{\displaystyle \frac{\partial}{\partial x^{i}}}.$$
In this case $A=0$ and $D=d_{E}-\delta$. 
If $u$ is in $^{E}T_{poly}(M)$ or $^{E}D_{poly}(M)$, a computation
shows that 
$$\lambda(u)= 
{\displaystyle \mathop \sum_{\alpha_{1}, \dots, \alpha_{d}}}
{\displaystyle \frac {(y^{1})^{\alpha_{1}}}{\alpha_{1}!}}\dots 
{\displaystyle \frac {(y^{d})^{\alpha_{d}}}{\alpha_{d}!}}
\left [\left (
{\displaystyle \frac{\partial}{\partial x^{1}}}\right )^{\alpha_{1}}
\dots 
{\displaystyle \left ( \frac{\partial}{\partial x^{d}}
\right )} ^{\alpha_{d}}\right ]\cdot u.$$
For example
$$\begin{array}{l}
\lambda_{T} \left (\gamma^{j_{1},\dots,j_{p}} 
{\displaystyle \frac{\partial}{\partial x^{j_{1}}}}\wedge \dots \wedge 
{\displaystyle \frac{\partial}{\partial x^{j_{p}}}} \right )\\
= 
{\displaystyle \mathop \sum_{\alpha_{1}, \dots, \alpha_{d}}}
{\displaystyle \frac {(y^{1})^{\alpha_{1}}}{\alpha_{1}!}}\dots 
{\displaystyle \frac {({y^{d})^{\alpha_{d}}}}{\alpha_{d}!}}
\left [
{\displaystyle \left (\frac{\partial}{\partial x^{1}}\right )^{\alpha_{1}}}
\dots 
{\displaystyle \left ( \frac{\partial }{\partial x^{d}} \right )^{\alpha_{d}}}
(\gamma^{j_{1}, \dots ,j_{p}})\right ]
{\displaystyle \frac{\partial}{\partial y^{j_{1}}}}\wedge \dots \wedge 
{\displaystyle \frac{\partial}{\partial y^{j_{p}}}}\\
= \lambda_{T} \left ( \gamma^{j_{1}, \dots ,j_{p}} \right )
{\displaystyle \frac{\partial}{\partial y^{j_{1}}}}\wedge \dots \wedge 
{\displaystyle \frac{\partial}{\partial y^{j_{p}}}}.
\end{array}$$

From lemma ~\ref{ouvert de coordonnees}, we see that 
$\overline{\cal V}^{B}=\overline{\cal V}$. 
If $a$ is in ${\cal O}_{\R^{d}}$, one has 
$${\displaystyle \frac{\partial}{\partial
    y^{i}}}\lambda(a)= 
\lambda \left ( {\displaystyle \frac{\partial a}{\partial  x^{i}}}
\right ).$$
Then it is easy to see that in this special case 
$\overline{\cal V} \circ \lambda_{T}$ 
takes its values in ${\cal D}_{poly} \bigcap Ker D$.

$B_{\Gamma}$ makes sense if we change the last argument by a
polydifferential operator with coefficients in $M$ and it is not hard
to see that 
$${\cal V}_{M}^{[n]}={\displaystyle \mathop \sum_{\Gamma \in G_{n+1,m}}}
W_{\Gamma}B_{\Gamma}.$$

\section{Applications}

In this section, we set $O=\Gamma ({\cal O}_{X})$. 
Let $E$ be a Lie algebroid,  ${\cal M}$ a $D(E)$-module and 
$M=\Gamma ({\cal M})$.
We denote by ${\cal V} _{\cal M}$ the quasi-isomorphism  of 
$L_{\infty}$-modules
over $\Gamma \left (^{E}T_{poly}^{*}\right )[[h]]$ from 
$\Gamma \left (^{E}T_{poly}^{*}({\cal M})\right )[[h]]$ 
to $\Gamma \left ( ^{E}D_{poly}^{*}({\cal M})\right )[[h]]$ 
given by theorem ~\ref{formality}. 
Then  ${\cal V}_{{\cal O}_{X}}={\Upsilon}$ is the 
$L_{\infty}$-quasi-isomorphism of DGLAs from 
$\Gamma \left ( ^{E}T_{poly}^{*}\right )[[h]]$ to 
$\Gamma \left ( ^{E}D_{poly}^{*}\right )[[h]]$ constructed by 
D. Calaque (~\cite{C1}).
 Let $\pi_{h}$ be a Maurer Cartan element of
$\Gamma \left (^{E}T_{poly}^{*}\right )[[h]]$. This means that 
$$\pi_{h}\in \Gamma \left (^{E}T_{poly}^{1}\right )[[h]] \;\;{\rm and}\;\; 
[\pi_{h}, \pi_{h}]_{S}=0. $$ Then it is well known that 
${\displaystyle \mathop \sum_{n \geq 1} \frac{1}{n!}
{\Upsilon} ^{[n]}(\pi_{h}, \dots, \pi_{h})}$ 
is a Maurer Cartan element of 
$\Gamma \left ( ^{E}D_{poly}^{*}\right )[[h]]$ 
(see ~\cite{AMM} p. 63). 
We set 
$$\Pi _{h}=1\otimes 1 + 
{\displaystyle \mathop \sum_{n \geq 1} \frac{1}{n!}
{\Upsilon } ^{[n]}(\pi_{h}, \dots, \pi_{h})}.$$
As $\Gamma \left (^{E}T^{*}_{poly}({\cal M}) \right )[[h]]$ 
is a module over the DGLA 
$\Gamma \left (^{E}T_{poly}^{*}\right )[[h]]$, the map 
$$\begin{array}{rcl}
\pi_{h} {\displaystyle \mathop \cdot_{S}}- : 
\Gamma \left (^{E}T_{poly}^{k}({\cal M})\right )[[h]] 
& \to & \Gamma \left ( ^{E}T_{poly}^{k+1}({\cal M})\right )[[h]]\\
y& \mapsto & \pi_{h} \cdot_{S}y
\end{array}$$
is a differential over 
$\Gamma \left (^{E}T^{*}_{poly}({\cal M}) \right )[[h]]$ 
(see ~\cite{D2} proposition 3 of section 2.3). Similarly,    
 $\Pi _{h} \cdot_{G}-$ defines a differential on 
$\Gamma \left (^{E}D_{poly}^{*}({\cal M})\right )[[h]]$.

\begin{proposition}
The map
$$\begin{array}{rcl}
\left ({\cal V}_{\cal M}\right )'_{\pi}:\;
\left (\Gamma \left ( ^{E}T_{poly}^{*}({\cal M})\right )[[h]], 
\pi_{h}\cdot_{S} - \right ) & \to & 
\left ( \Gamma \left (^{E}D_{poly}^{*}({\cal M})\right )[[h]], 
\Pi_{h}\cdot_{G} - \right )\\
y & \mapsto & 
{\displaystyle \mathop \sum _{p=0}^{\infty}
\frac{1}{p!}
{\cal V}_{{\cal M}}^{[p]}(\pi_{h}, \dots  ,\pi_{h}, y )}
\end{array}$$
is a quasi-isomorphism.
\end{proposition}
{\it Proof of the proposition : }

The proposition follows from proposition 3 of paragraph 2.3
of ~\cite{D2} and the definition of the 
$L_{\infty}$-module structure of 
$\Gamma \left (^{E}D^{*}_{poly}({\cal M})\right )$
 over $\Gamma \left (^{E}T^{*}_{poly}\right )$.$\Box$\\

If $E$ is a Lie algebroid equipped with an $E$-bivector 
$\pi \in \Gamma ( \bigwedge ^{2} E)$ satisfying $[\pi, \pi]=0$, it is
called a Poisson Lie algebroid. If $E=TX$, we recover Poisson manifolds. 
Then, one can construct a Lie algebroid structure on $E^{*}$ 
in the
following way. Let $\pi^{\sharp}$ be the bundle map from $E^{*}$ to $E$
associated to $\pi$ and 
$\omega_{*}=\omega \circ \pi^{\sharp} : E^{*} \to TX$. 
Define a Lie bracket on
$E^{*}$ by :
$$\forall \theta , \eta \in E ^{*}, \;\;
[\theta , \eta ]=L_{\pi^{\sharp}\theta}(\eta) - 
L_{\pi^{\sharp}\eta}(\theta) - d \pi (\theta , \eta)$$
where $L$ denotes the Lie derivative.
Then $E^{*}$, endowed with the bracket above and the anchor 
$\omega_{*}$ is a Lie algebroid (~\cite{KM}, ~\cite{MX}) 
and $E$ is a Lie bialgebroid. 
The differential of the Lie cohomology complex of $E^{*}$ is 
$d_{*}=[\pi, - ] : \Gamma (\bigwedge ^{k} E)\to 
\Gamma (\bigwedge^{k+1}E)$.

Assume that we are in the case where $E$ is a Poisson Lie algebroid with
Poisson bivector $\pi$. Then, in the proposition above one may take 
$\pi_{h}=h \pi$ and 
Calaque (~\cite{C1}) shows that $\Pi_{h}$ is a twistor for the standard Hopf 
algebroid $U(\Gamma (E))[[h]]$ (~\cite{X}). 

From now on, we assume that $E=TX$ and that 
$\pi$ is a Poisson bracket on $X$. Then
the twistor $\Pi_{h}$ defines a star product on $O[[h]]$ (~\cite{X})
in the following way
$$\forall (f, g) \in O ,\;\; \Pi_{h}(f,g)=f*_{h}g.$$
Set 
$$f*_{h}g=fg + 
{\displaystyle \mathop \sum_{i=1}^{\infty}}B_{i}(f,g)h^{i}$$

\begin{proposition}
$M[[h]]$ can be
endowed with an $O[[h]]\otimes O[[h]]^{op}$-module
structure as follows : for all $a$ in $O$ and all $m$ in $M$, 
$$a * m = a\cdot m + 
{\displaystyle \mathop \sum_{i=1}^{\infty}h^{i}B_{i}(a,-  )}\cdot m,
\;\;\;
m*a = a \cdot m + 
{\displaystyle \mathop \sum_{i=1}^{\infty}}h^{i}B_{i}(-, a)\cdot m $$
\end{proposition}
{\it Proof of the proposition :}

The proof of the proposition is a straightforward verification using
the associativity of the star product. $\Box$\\

Applying the exact functor $N \mapsto N[[h]]$, we get an injection 
$$\Gamma \left (^{E}D_{poly}^{k}({\cal M})\right )[[h]] \hookrightarrow 
Hom_{{\R} [[h]]} 
\left (O[[h]]^{\otimes_{\R [[h]]}^{k+1}} , M[[h]] \right
).$$
The image of $\Gamma \left ( ^{E}D_{poly}^{*}({\cal M})\right )[[h]]$ in 
$Hom_{\R [[h]]} 
\left (O[[h]]^{\otimes_{\R [[h]]}^{*+1}} , M[[h]] \right
)$ will be denoted 
$Homdiff_{\R [[h]]} 
\left (O[[h]]^{\otimes_{\R [[h]]}^{*+1}} , M[[h]] \right )$ .

Recall that the Hochschild cohomology of $O[[h]]$ with values in the
bimodule $M[[h]]$, $HH^{*} \left (O[[h]], M[[h]] \right )$, is the
cohomology of the complex \linebreak 
$\left ( Hom_{\R [[h]]}(O[[h]]^{\otimes^{*}_{\R [[h]]}}, 
M[[h]]), \beta \right )$
where, with obvious notations, 
$$\begin{array}{rcl}
\beta(\lambda)(a_{1}, \dots, a_{n+1})& =& 
a_{1}*\lambda(a_{2}, \dots,a_{n+1})\\
&+&{\displaystyle \mathop \sum_{0<i<n+1}}(-1)^{i}
\lambda (a_{1}, \dots, a_{i}*a_{i+1},\dots, a_{n+1})\\
&+& (-1)^{n+1}\lambda(a_{1}, \dots, a_{n})*a_{n+1}.\\
\end{array}$$
Denote by $HH_{md}^{*}\left (O[[h]],M[[h]] \right )$ the cohomology of the 
complex \linebreak
$\left ( Homdiff_{\R[[h]]}(O[[h]]^{\otimes^{*}}, M[[h]]), \beta\right )$.

The complex 
$\left ( \Gamma ( ^{E}T_{poly}^{*}({\cal M}))[[h]], \pi_{h}\cdot_{S}
\right )$
computes 
the Lichnerowicz-Poisson cohomology of the $\R [[h]]$-Poisson algebra 
(defined by the bivector $\pi_{h}$) 
$O[[h]]$ with coefficients in $M[[h]]$, 
$H^{i}_{Poisson} \left ( O[[h]], M[[h]]\right )$(~\cite{Li}, ~\cite{Hu}). 
Denote by \linebreak
$H^{i}_{Poisson} \left ( O[[h]], M[[h]]\right )$ the
Lichnerowicz-Poisson cohomology of the $\R [[h]]$-Poisson algebra 
(defined by the bivector $\pi_{h}$) 
$O[[h]]$ with coefficients in $M[[h]]$ (~\cite{Hu}). 
It is computed by the complex 
$\left ( \Gamma ( ^{E}T_{poly}^{*}({\cal M}))[[h]], \pi_{h}\cdot_{S} \right )$.
 The complex \linebreak  
$\left ( \Gamma (^{E}D_{poly}^{*}({\cal M}))[[h]], \Pi_{h}
  \cdot_{G}\right )$ 
computes 
$HH_{md}^{*}\left ( O[[h]], M[[h]] \right )$. 
We get the following corollary :

\begin{corollary}
One has an isomorphism 
$$H^{i}_{Poisson}\left ( O[[h]], M[[h]]\right )\simeq 
HH^{i}_{md}\left ( O[[h]], M[[h]] \right ).$$
\end{corollary}  

The exterior product, which will be denoted by $\wedge$, endows 
$H^{*}\left ( \Gamma \left ( ^{E}T_{poly}^{*}\right ), [\pi_{h}, \cdot
]\right )$ with an associative supercommutative algebra structure. It also
endows \linebreak
$H^{*}\left ( \Gamma \left ( ^{E}T_{poly}^{*}({\cal M})\right ), 
\pi_{h} \cdot_{S} \right )$ with a 
$\left [ H^{*}\left ( \Gamma \left ( ^{E}T_{poly}^{*}\right ), 
[\pi_{h}, \cdot ]\right ), \wedge \right ]$-module structure. 

To simplify the notation, from now on, we write $\Pi$ instead of $\Pi_{h}$.  
$D^{*}_{poly}$ is endowed with an associative graded product,
$\sqcup_{\Pi}$, compatible with the differential $[\Pi, \cdot ]$ 
defined by : 

$$\begin{array}{l}
\forall t_{1} \in \Gamma \left ( D_{poly}^{k_{1}-1} \right ), \; 
\forall t_{2} \in \Gamma \left ( D_{poly}^{k_{2}-1} \right ),\; 
\forall a_{1}, \dots, a_{k_{1}+k_{2}} \in O,\\
(t_{1}\sqcup_{\Pi}t_{2})(a_{1}, \dots, a_{k_{1}+k_{2}})=
t_{1}(a_{1}, \dots, a_{k_{1}}) \star_{h} t_{2}(a_{k_{1}+1}, \dots, 
a_{k_{1}+k_{2}})\\
\end{array}$$ 
Thus, 
$\left [ H^{*}\left (\Gamma (D_{poly}), [\Pi, \cdot] \right ), \sqcup_{\Pi}
\right ]$ is an associative graded algebra.
$t_{1}\sqcup_{\Pi} t_{2}$ is also  \linebreak
defined if 
$t_{2} \in \Gamma \left ( D_{poly}^{k_{2}-1}({\cal M}) \right )$. Thus,    
$\left [ H^{*}\left (\Gamma (D_{poly}({\cal M})), \Pi \cdot_{G} \right ), 
\sqcup_{\Pi} \right ]$ is a \linebreak 
$\left [ H^{*}\left (\Gamma (D_{poly}), [\Pi, \cdot] \right ), \sqcup_{\Pi}
\right ]$-module.  \\

If $X=\R^{d}$ and $E=T\R^{d}$, Kontsevich has proved 
(~\cite{Ko}, see ~\cite{MT} for a detailed proof) 
that the algebras 
$\left [ H^{*}\left ( \Gamma \left ( T_{poly}^{*}\right ), 
[\pi_{h}, \cdot ]\right ), \wedge \right ]$ and 
$\left [ H^{*}\left (\Gamma (D_{poly}), [\Pi, \cdot] \right ), \sqcup_{\Pi}
\right ]$ are isomorphic. We will extend this result to our case. \\

{\bf Remark :}
In ~\cite{CFT}, 
a star product $*$ is constructed on any manifold $X$ so that 
the algebras 
$\left  [ H^{0}\left ( \Gamma \left ( T_{poly}^{*}\right ), 
[\pi_{h}, \cdot ]\right ), \wedge \right ]$ and 
$\left [ H^{0}\left (\Gamma \left ( D_{poly} \right ), 
[*, \cdot] \right ), \sqcup_{\Pi}
\right ]$ are isomorphic.\\

\begin{theorem}\label{compatibility}
Assume that $X=\R^{d}$ and 
$E=T\R^{d}$. The \linebreak 
$\left [ H^{*}\left ( \Gamma \left ( T_{poly}^{*}\right ), 
[\pi_{h}, \cdot ]\right ), \wedge \right ] \simeq 
\left [ H^{*}\left (\Gamma (D_{poly}), [\Pi, \cdot] \right ), \sqcup_{\Pi}
\right ]$-modules \linebreak 
$\left [ H^{*}\left ( \Gamma \left ( T_{poly}^{*}({\cal M})\right ), 
\pi_{h} \cdot_{S}\right ), \wedge \right ]$ and 
$\left [ H^{*}\left (\Gamma \left ( D_{poly}({\cal M})\right ), 
\Pi \cdot_{G} \right ), \sqcup_{\Pi}
\right ]$ are isomorphic.
\end{theorem} 

{\it Proof of the theorem \ref{compatibility}}

In this proof, we keep the notations of the proof of the formality
theorem (paragraph 4.3).  
We could reproduce the proof of ~\cite{MT} using the explicit expression we
found for ${\cal V}_{\cal M}$ in the paragraph 4.4. We will use the
decomposition 
${\cal V}_{\cal M}=\lambda_{D}^{-1}\circ \overline{\cal V} \circ \lambda_{T}$ 
and use the results of ~\cite{MT}. 
Put 
$$\overline{\Pi}={\displaystyle \mathop \sum_{n\geq 1}}
{\cal U}^{[n]}
\left ( \lambda_{T}(\pi_{h}), \dots, \lambda_{T}(\pi_{h}) \right ).$$ 
\begin{lemma}\label{MT}
Let $k_{1}$ and $k_{2}$ be in $\N$. 
If $\tau_{1} \in \Gamma \left ({\cal T}_{poly}^{k_{1}-1}\right )$,  
$\tau_{2} \in \Gamma \left ( {\cal T}_{poly}^{k_{2}-1}({\cal M})\right )$ 
and $m=k_{1}+k_{2}$ then one
has 
 $$\begin{array}{l}
\overline{\cal V}'_{\lambda_{T}(\pi_{h})}(\tau_{1} \wedge \tau_{2}) -
\overline{\cal U}'_{\lambda_{T}(\pi_{h})}(\tau_{1})
\sqcup_{\overline{\Pi}}
\overline{\cal V}'_{\lambda_{T}(\pi_{h})}(\tau_{2})=\\
{\displaystyle \sum_{n \geq 0}}{\displaystyle \frac{h^{n}}{n!}}
{\displaystyle \sum_{\Delta \in G_{n+2, m-1}}}a_{\Delta}
\overline{\Pi}\cdot_{G}B_{\Delta}
(\lambda_{T} (\pi)\otimes \dots \otimes \lambda_{T}(\pi)\otimes 
\tau_{1}\otimes \tau_{2}) \\
+{\displaystyle \sum_{n \geq 0}}{\displaystyle \frac{h^{n}}{n!}}
{\displaystyle \sum_{\Delta \in G_{n+1, m}}}b_{\Delta}
 (-1)^{(k_{1}-1)k_{2}}
B_{\Delta}(\lambda_{T}(\pi)\otimes \dots \otimes \lambda_{T}(\pi)
\otimes [\lambda_{T}(\pi), \tau_{1}] \otimes \tau_{2}) \\
{\displaystyle \sum_{n \geq 0}}{\displaystyle \frac{h^{n}}{n!}}
{\displaystyle \sum_{\Delta \in G_{n+1, m}}}b_{\Delta}
(-1)^{k_{1}(k_{2}-1)}
B_{\Delta}\left ( \lambda_{T}(\pi)\otimes \dots  \lambda_{T}(\pi)\otimes
  \tau_{1}\otimes  \lambda_{T}(\pi)\cdot_{S}\tau_{2} \right )
\end{array}$$
where $a_{\Delta}$ and $b_{\Delta}$ are real. 
\end{lemma}

{\it Proof of the lemma \ref{MT} :}

The lemma \ref{MT} is proved for ${\cal M}={\cal O}_{X}$ in
~\cite{MT}. 
 Actually, the formula of lemma ~\ref{MT} is slightly different from that of 
~\cite{MT}. To get it, one has to reproduce the proof of ~\cite{MT} and
make play to the vertices $n-1$ and $n$ the role played by the
vertices $1$ and $2$.
Hence the lemma ~\ref{MT}  holds for 
$\tau_{2}$ in 
$\Gamma \left ( {\cal T}^{k_{2}-1}_{poly}\right )
{\displaystyle \mathop \otimes_{O}}M$. 
We will now
prove that it is true for $\tau_{2}$ in 
$\Gamma \left ({\cal T}_{poly}^{k_{2}-1}({\cal M})\right )$. 
If we apply it to $(f_{1}, \dots, f_{m})$ in 
$\R[[y^{1},\dots, y^{d}]]^{m}$, 
the relation of the lemma \ref{MT} can be written 
${\displaystyle \sum_{n \geq 0}}h^{n}F_{n}= 
{\displaystyle \sum_{n \geq 0}}h^{n}G_{n}$   where the 
$F_{n}$'s and the $G_{n}$'s are maps from 
$\Gamma \left ( {\cal T}^{k_{2}-1}_{poly}\right )
{\displaystyle \mathop \otimes_{O}M}$ 
to $M[[y^{1},\dots,y^{d}]]$. 
Let $I$ be the ideal of $O[[y^{1}, \dots, y^{d}]]$
generated by $y^{1}, \dots, y^{d}$. The $F_{n}$'s and the $G_{n}$'s
are continuous for the $I$-adic topology. This is a consequence of  the
following two remarks.

${\bullet}$ Let $\gamma_{1}, \dots, \gamma_{p}$ be elements of 
$\Gamma \left ({\cal T}_{poly}\right )$ 
and let $(g_{1},\dots, g_{m})$ be elements of 
$O[[y^{1},\dots, y^{d}]]$. 
Let $\Gamma$ be an admissible graph of type
$(p+1,m)$. The map 
$$\begin{array}{rcl}
\Gamma \left ( {\cal T}^{k_{2}-1}_{poly}\right )
{\displaystyle \mathop \otimes_{O}M}
& \to & 
M[[y^{1}, \dots, y^{d}]]\\
\mu &\mapsto& B_{\Gamma}(\gamma_{1}, \dots, \gamma_{p}, \mu)
(g_{1},\dots, g_{m})
\end{array}$$
is continuous for the $I$-adic topology as it sends 
$I ^{N}\Gamma \left ( {\cal T}^{k_{2}-1}_{poly}\right )
{\displaystyle \mathop \otimes_{O}M}$ to \linebreak 
$I^{N-p}M[[y^{1},\dots,y^{d}]]$. 

$\bullet$ Let $\Gamma$ be an admissible graph of type
$(p,2)$ and let $g$ be an element of 
$O[[y^{1}, \dots, y^{d}]]$. The map 
$$\begin{array}{rcl}
O[[y^{1}, \dots, y^{d}]]
{\displaystyle \mathop \otimes_{O}M}& \to & 
M[[y^{1},\dots, y^{d}]]\\
\mu &\mapsto& B_{\Gamma}(\lambda_{T} (\pi), \dots, \lambda_{T}(\pi))(f, \mu)
\end{array}$$
is continuous for the $I$-adic topology as it sends 
$I^{N}O[[y^{1}, \dots, y^{d}]]
{\displaystyle \mathop \otimes_{O}M}$ to \linebreak
$I^{N-p}M[[y^{1},\dots,y^{d}]]$. 

This finishes the proof of the lemma \ref{MT}.$\Box$ \\

Now, we go back to the
proof of the theorem \ref{compatibility}

Let $t_{1}$ be in $\Gamma (T^{k_{1}-1}_{poly})[[h]]\cap Ker [\pi_{h},]$ and 
$t_{2}$ be in 
$\Gamma \left (T^{k_{2}-1}_{poly}({\cal M})\right )[[h]]
\cap Ker (\pi_{h} \cdot_{S})$. 
We apply the lemma ~\ref{MT} to $\tau_{1}=\lambda_{T} (t_{1})$ and 
$\tau_{2}=\lambda_{T}^{\cal M} (t_{2})$. We get 
 
$$\begin{array}{l}
\overline{\cal V}'_{\lambda_{T}(\pi_{h})}
(\lambda_{T}(t_{1}) \wedge \lambda_{T}^{\cal M}(t_{2})) -
\overline{\cal U}'_{\lambda_{T}(\pi_{h})}(\lambda_{T}(t_{1}))
\sqcup_{\overline{\Pi}}
\overline{\cal V}'_{\lambda_{T}(\pi_{h})}(\lambda_{T}^{\cal M}(t_{2}))\\
={\displaystyle \sum_{n \geq 0}}{\displaystyle \frac{h^{n}}{n!}}
{\displaystyle \sum_{\Delta \in G_{n+2, m-1}}}a_{\Delta}
\overline{\Pi}\cdot_{G}B_{\Delta}
(\lambda_{T}(\pi)\otimes \dots \otimes \lambda_{T}(\pi)\otimes 
\lambda_{T}(t_{1})\otimes \lambda_{T}^{\cal M}(t_{2})). 
\end{array}$$
Apply $(\lambda_{D}^{\cal M})^{-1}$ and use the following facts :

$\bullet$ $\lambda_{D}^{-1}(\overline{\Pi})=\Pi$ .

$\bullet $ With obvious notations, one has. 
$$\lambda_{D}(\sigma_{1})  
\sqcup_{\overline{\Pi}}\lambda _{D}^{\cal M}(\sigma_{2})
=\lambda_{D}^{\cal M}(\sigma_{1}  \sqcup_{\Pi} \sigma_{2}).$$

$\bullet$ $B_{\Delta}\left ( \lambda_{T}(\pi), \dots,
  \lambda_{T}(\pi), \lambda_{T}(t_{1}), \lambda_{T}^{\cal M}(t_{2}) \right )=
\lambda_{D}^{\cal M}\left ( B_{\Delta}\left ( \pi, \dots, \pi, t_{1}, t_{2}
  \right ) \right )$

We get 

$$\left ({\cal V}_{\cal M}\right )'_{\pi}(t_{1} \wedge t_{2}) -
{\cal U}'_{\pi}(t_{1})
\sqcup_{\Pi}\left ({\cal V}_{\cal M}\right )'_{\pi}(t_{2})=\\
{\displaystyle \sum_{n \geq 0}}{\displaystyle \frac{h^{n}}{n!}}
{\displaystyle \sum_{\Delta \in G_{n+2, m-1}}}a_{\Delta}
 \Pi\cdot_{G}B_{\Delta}
(\pi\otimes \dots \otimes \pi\otimes 
t_{1}\otimes t_{2}). $$
The right hand side is a coboundary for the Hochschild cohomology
complex. This finishes the proof of the theorem \ref{compatibility}.$\Box$\\

{\bf Remark :}

Assume that $X$ is the dual of a real Lie algebra endowed with its
Kirillov-Kostant-Souriau Poisson structure denoted by $\pi$. 
Recall that if $\xi$ and $\eta$ are elements of $\goth{g }$
considered as linear forms  on $\goth{g} ^{*}$, then  
$$\pi (\xi , \eta) =[\xi, \eta ].$$
If $M={\cal O}_{X}$, 
the isomorphism given by theorem ~\ref{compatibility}
has been studied. 
If $i=0$, it gives Duflo's isomorphism (~\cite{Du}, ~\cite{Ko}). 
By analyzing which graphs contributes to 
$\left ( {\cal V}_{\cal M} \right )_{\pi}'$, 
Pevsner and Torossian (~\cite{PT}) have shown that  
that Duflo's isomorphism
extends to an isomorphism from 
$H^{*}_{Poisson}({\goth g}, S({\goth g}))$ to 
$H^{*}({\goth g}, U({\goth g}))$.\\

\noindent Sophie Chemla\\
UPMC Universit{\'e} Paris 6\\
UMR 7586\\
Institut de math{\'e}matiques\\
75005 Paris\\
France.\\
schemla@math.jussieu.fr


\begin{thebibliography}{BCKT}



\bibitem[AK]{AK} R. Almeida, A. Kumpera, {\it Structure produit dans la
  cat{\'e}gorie des alg{\'e}bro{\"\i}des de Lie}, Ann. Acad. Brasil. Cienc
{\bf 53} (1981), 247--250. 

\bibitem[AMM]{AMM} D. Arnal, D. Manchon and M. Masmoudi,   
{\it Choix des signes pour la  formalit{\'e} de M. Kontsevich}, 
Pacific Journal of math {\bf 203}, 1, (2002).

\bibitem[Bo]{Bo} A. Borel : {\it Algebraic D-modules}, Academic press, 1987. 

\bibitem[BCKT]{BCKT} A. Brugui{\`e}res, A. Cattaneo and  B. Keller and 
C. Torossian, {\it D{\'e}formation, Quantification, Th{\'e}orie de Lie}, 
Panoramas et synth{\`e}se, SMF (2005). 

\bibitem [C1]{C1} D. Calaque, {\it Formality for Lie algebroids}, 
Communications in mathematical physics, {\bf 257},3 (2005), 563--578. 

\bibitem [C2]{C2} D.Calaque, {\it Th{\'e}or{\`e}mes de formalit{\'e} pour les
  alg{\'e}bro{\"\i}des  de Lie et quantification des $r$-matrices
  dynamiques}, Th{\`e}se de l'IRMA.

\bibitem [CDH]{CDH} D. Calaque,V. Dolgushev and G. Halbout,
 {\it Formality theorem
  for Hochschild chains in the Lie algebroid setting}, 
J. Reine Angew. Math. {\bf 612} (2007), 81--127 (math KT/0504372).

\bibitem[CFT]{CFT} A.S. Cattaneo, G. Felder and L.Tomassini,  
{\it From local to global deformation quantization of Poisson
  manifolds}, 
Duke Math. J. {\bf 115}, 2 (2002), 329--352. 

\bibitem[Ch1]{Ch1} S. Chemla,  
{\it Poincar{\'e} duality for $k-A$-Lie superalgebras}, 
Bulletin de la Soci{\'e}t{\'e} math{\'e}matique de France, {\bf 122} (1994),
371--397 .

\bibitem[Ch2]{Ch2} S. Chemla, 
{\it A duality property for complex Lie algebroids},
Mathematische Zeitscrift {\bf 232} (1999), 367-388.  

\bibitem[Ch3]{Ch3} S. Chemla,
 {\it An inverse image functor for Lie algebroids},
Journal of algebra {\bf 269} (2003), 109--135.


\bibitem[D1]{D1}  V. Dolgushev,  
{\it Covariant and equivariant formality theorem}, 
Adv. Math. {\bf 191}, 1 (2005) 147--177 (math.QA/0307212). 

\bibitem[D2]{D2}  V. Dolgushev, 
{\it A formality theorem for Hochschild chains},
Adv. Math. {\bf 200}, 1 (2006), 51--101 (math QA/0402248).

\bibitem[D3]{D3} V. Dolgushev, 
{\it A proof of Tsygan's formality conjecture for arbitrary smooth manifolds},
PhD Thesis, math QA/0504420.

\bibitem[Du]{Du} M. Duflo, {\it Op{\'e}rateurs diff{\'e}rentiels bi-invariants sur
  un groupe de Lie}, Annales Scientifiques de l' Ecole  Normale 
Sup{\'e}rieure {\bf 10},  (1977), 107--144.

\bibitem[ELW]{ELW} S. Evens, J-H Lu, A. Weinstein,
{\it Transverse measures, the
modular class and a cohomology pairing for Lie algebroid}, Quarterly
Journal of Math, {\bf 50} (1999), 171--220.

\bibitem[Fe]{Fe} B. Fedosov, 
{\it A simple geometrical construction of deformation quantization}, 
J. Diff. Geom. {\bf 40} (1994) 213-238.

\bibitem[F]{F} R.L. Fernandes, 
{\it Lie algebroids, holonomy and characteristic
  class}, Adv. Math. {\bf 170}, 1 (2002), 119--179.

\bibitem[HKR]{HKR} G. Hochschild, B. Kostant and A. Rosenberg,
 {\it Differential forms on regular affine algebras }, 
Trans. Amer. Math. Soc. {\bf 102} (1962), 383--408.

\bibitem [HS]{HS} V. Hinich and V. Schechtman, 
{\it Homotopy Lie algebras}, I.M
Gelfand seminar, Adv. Sov. Math. {\bf 16}, 2 (1993), 1--28.  

\bibitem[Hu]{Hu} J. Huebschmann, {\it Poisson cohomology and quantization},
J. reine angew. Math. {\bf 408} (1990), 57--113.

\bibitem[Ka]{Ka} M. Kashiwara, {\it D-modules and microlocal
    calculus}, Translations of mathematical monographs, {\bf 217}, 
American Mathematical Society (2003). 

\bibitem[Ko]{Ko} M. Kontsevitch,  {\it Deformation quantization of Poisson
  manifolds}, Lett. Math. Phys. {\bf 66}, 3 (2003), 157--216.

\bibitem[KM]{KM} Y. Kosmann-Schwarzbach and F. Magri, {\it Poisson Nijenhuis
  structures}, Annales de l'Institut Henri Poincar{\'e}, S{\'e}rie A, {\bf
53} (1990), 35-81.

\bibitem[KS]{KS} M. Kashiwara- P. Schapira, {\it Sheaves on manifolds}, 
Grundlehren  der mathematischen Wissenschaften, A series of
Comprehensive Studies in Mathematics, Springer-Verlag  (1994). 

\bibitem[LS]{LS} T. Lada and J. Stasheff, {\it Introduction to SH Lie algebras
  for physicists}, Intern. J. Theor. Phys. {\bf 32} (1993), 1087--1103.

\bibitem[Li]{Li} A. Lichnerowicz, {\it Les vari{\'e}tes de Poissons et leurs
  alg{\`e}bres associ{\'e}es}, Journal of differential geometry {\bf 12}
(1977), 253--300.

\bibitem[MT]{MT} D. Manchon and C. Torossian, {\it Cohomologie tangente et
  cup-produit pour la quantification de Kontsevitch}, Ann. Math. Blaise
  Pascal {\bf 10}, 1 (2003), 75-106.

\bibitem [MX]{MX} K.C.H Mackenzie and Ping Xu,  
{\it Lie bialgebroids and Poisson  groupoids}, 
Duke Math. Journal {\bf 73} (1994), 415--452.

\bibitem[PT]{PT} M. Pevsner and C. Torossian, {\it Isomorphisme de Duflo et
  cohomologie tangentielle}, Journal of geometry and Physics, 
  {\bf 51}, 4 (2004), 486--505.

\bibitem[R]{R} G.S.Rinehart, {\it Differential form on general
    commutative algebra}, Trans. Amer. Math. Soc {\bf 108} (1963),
195--222. 

\bibitem[V]{V} J. Vey, {\it D{\'e}formation du crochet de Poisson sur une
  vari{\'e}t{\'e} symplectique}, Comment. Math. Helv. {\bf 50} (1975), 
421--454. 

\bibitem [X]{X} P. Xu,  {\it Quantum groupoids}, 
Comm. Math. Phys. {\bf 206} (2001), 
539--581. \\

\end{thebibliography}
\end{document}